\documentclass[12pt, reqno]{amsart}
\usepackage{amssymb, amstext, amscd, amsmath}

\usepackage{xy}
\xyoption{all}

\makeatletter
\def\@cite#1#2{{\m@th\upshape\bfseries%
[{#1\if@tempswa{\m@th\upshape\mdseries, #2}\fi}]}} \makeatother
\theoremstyle{plain}
\newtheorem{thm}[subsection]{Theorem}
\newtheorem{cor}[subsection]{Corollary}
\newtheorem{prop}[subsection]{Proposition}
\newtheorem{lem}[subsection]{Lemma}
\theoremstyle{definition}
\newtheorem{rem}[subsection]{Remark}
\newtheorem{example}[subsection]{Example}

\newtheorem{defn}[subsection]{Definition}

\newtheorem{eg}[subsection]{Example}

\newcommand{\bC}{{\mathbb{C}}}

\newcommand{\bF}{{\mathbb{F}}}

\newcommand{\bN}{{\mathbb{N}}}

\newcommand{\bT}{{\mathbb{T}}}

\newcommand{\bZ}{{\mathbb{Z}}}

  \newcommand{\A}{{\mathcal{A}}}
  \newcommand{\B}{{\mathcal{B}}}
  \newcommand{\C}{{\mathcal{C}}}

  \newcommand{\F}{{\mathcal{F}}}
  
\renewcommand{\H}{{\mathcal{H}}}

  \newcommand{\K}{{\mathcal{K}}}
\renewcommand{\L}{{\mathcal{L}}}
\newcommand{\M}{{\mathcal{M}}}

\renewcommand{\S}{{\mathcal{S}}}

  \newcommand{\Z}{{\mathcal{Z}}}

\newcommand{\fA}{{\mathfrak{A}}}

\newcommand{\fF}{{\mathfrak{F}}}
\newcommand{\fG}{{\mathfrak{G}}}

\newcommand{\fJ}{{\mathfrak{J}}}

\newcommand{\fg}{{\mathfrak{g}}}

\newcommand{\fs}{{\mathfrak{s}}}
\newcommand{\ft}{{\mathfrak{t}}}
\newcommand{\fu}{{\mathfrak{u}}}

\newcommand{\Bk}{{\mathbf{k}}}
\newcommand{\Bl}{{\mathbf{l}}}
\newcommand{\Bm}{{\mathbf{m}}}
\newcommand{\Bn}{{\mathbf{n}}}

\newcommand{\rC}{{\mathrm{C}}}

\newcommand{\ep}{\varepsilon}
\renewcommand{\phi}{\varphi}
\newcommand{\upchi}{{\raise.35ex\hbox{\ensuremath{\chi}}}}

\newcommand{\qand}{\quad\text{and}\quad}

\newcommand{\qfor}{\quad\text{for}\quad}
\newcommand{\qforal}{\quad\text{for all}\quad}

\newcommand{\AND}{\text{ and }}

\newcommand{\dcup}{\operatorname{\,\dot\cup\,}}
\newcommand{\diag}{\operatorname{diag}}

\newcommand{\id}{{\operatorname{id}}}

\newcommand{\ran}{\operatorname{Ran}}
\newcommand{\spn}{\operatorname{span}}

\newcommand{\Ath}{\A_\theta}

\newcommand{\bsl}{\setminus}
\newcommand{\ca}{\mathrm{C}^*}
\newcommand{\cenv}{\mathrm{C}^*_{\text{env}}}

\newcommand{\Fth}{\mathbb{F}_\theta^+}
\newcommand{\Ftheta}{\mathbb{F}_\theta^+}

\newcommand{\Fockth}{\ell^2(\Fth)}

\newcommand{\mt}{\varnothing}

\newcommand{\ol}{\overline}
\newcommand{\one}{{\boldsymbol{1}}}

\newcommand{\wh}{\widehat}

\newcommand{\ltwo}{\ell^2}

\begin{document}
\title[Higher Rank Graph Algebras]
{Representations of \\Higher Rank Graph Algebras}

\author[K.R. Davidson {\protect \and} D. Yang]
{Kenneth R. Davidson {\protect \and} Dilian Yang}
\address{Kenneth R. Davidson, 
Pure Mathematics Department,
University of Waterloo, 
Waterloo, ON N2L 3G1,
CANADA}
\email{krdavids@uwaterloo.ca}
\address{Dilian Yang,
Pure Mathematics Department,
University of Waterloo, 
Waterloo, ON N2L 3G1,
CANADA}
\email{dyang@math.uwaterloo.ca}

\begin{abstract}
Let $\Fth$ be a $\Bk$-graph on a single vertex.
We show that every irreducible atomic $*$-representation is the 
minimal $*$-dilation of a group construction representation.
It follows that every atomic representation decomposes as a direct sum or integral
of such representations. We characterize periodicity of $\Fth$
and identify a symmetry subgroup $H_\theta$ of $\bZ^\Bk$.
If this has rank $s$, then $\ca(\Fth) \cong \rC(\bT^s) \otimes \fA$
for some simple C*-algebra $\fA$.
\end{abstract}

\thanks{First author partially supported by an NSERC grant.\\ \indent
Second author partially supported by the Fields Institute.}

\date{}
\maketitle

\section{Introduction}

There has been a lot of recent interest in the structure of operator algebras
associated to graphs (see \cite{Raeburn}).  
Kumjian and Pask \cite{KumPask} have introduced
the notion of higher rank graphs, which have a much more involved
combinatorial structure.  
The C*-algebras of higher rank graphs are widely studied
\cite{FMY,FowSim,KumPask2,PRRS,RaeSimYee1,RaeSimYee2,RS,SZ}.
Kribs and Power \cite{KP1} initiated the study of the
associated nonself-adjoint algebras.  
Power \cite{Power1} began a detailed study
of these operator algebras associated to higher rank
graphs with a single vertex.  This effort was continued with the authors
of this paper \cite{DPY, DPYdilation, DYperiod} with a detailed analysis
of rank 2 graphs on a single vertex.  Two important accomplishments there
were a complete structure theory for the atomic $*$-representations of
the 2-graph, and a characterization of periodicity leading to the structure
of the 2-graph C*-algebra in the periodic case.

The purpose of this paper is to extend those results to the case of $\Bk$-graphs
on a single vertex.  These objects form an interesting class of semigroups
with cancellation and unique factorization.  The combinatorial structure
of a $\Bk$-graph $\Fth$ is much more difficult to classify for $\Bk \ge 3$, 
but we do show that there are lots of examples.  

The first main result concerns atomic $*$-representations.  
The $*$-representations are the row isometric representations of $\Fth$
which are defect free (see the next section for definitions).
These are the representations of the semigroup which yield
$*$-representations of the associated C*-algebra $\ca(\Fth)$.
An important class of such representations are obtained as
inductive limits of the left regular representation of $\Fth$.
These have the additional property that there is an orthonormal basis
which is permuted (up to scalars of modulus 1) by the isometries
(which are the images of elements of $\Fth$).  In general, a representation
with such a basis is called atomic.

The analysis of these representations relies on dilation theory.
Every defect free, row contractive representation dilates to a
unique minimal $*$-dilation \cite{DPYdilation,SZ}.
So one can understand a $*$-representation by understanding
its restriction to a coinvariant cyclic subspace.
The key to our analysis is to show that there is a natural family
of atomic defect free representations modelled on the representations
of an abelian group of rank $\Bk$.
In the rank 2 case \cite{DPY}, a detailed case by case analysis led to
the conclusion that every irreducible atomic $*$-representation is
the dilation of one of these group constructions.
In this paper, we provide a direct argument that avoids the case by
case approach.  So it sheds new light even when $\Bk=2$.

Every $\Bk$-graph $\Fth$ has a faithful $*$-representation of inductive type.
This representation has a natural symmetry group $H_\theta \le \bZ^\Bk$.
The graph is aperiodic if and only if $H_\theta = \{0\}$.
In general, this is a free group of rank $s \le \Bk$.
Building on the detailed structure of periodic 2-graphs in \cite{DYperiod},
we show that there is a simple C*-algebra $\fA$ so that
$\ca(\Fth) \cong \rC(\bT^s) \otimes \fA$.

\section{Background}\label{S:hrg}

Kumjian and Pask \cite{KumPask} define a $\Bk$-graph as a small category
$\Lambda$ with a degree map $\deg : \Lambda\to\bN^{\Bk}$ satisfying the
\textit{factorization property}: for every $\lambda\in\Lambda$
and $m,n\in\bN^{\Bk}$ with $\deg (\lambda)=m+n$, there are unique
elements $\mu,\nu\in\Lambda$ such that
$\lambda=\mu\nu$ and  $\deg(\mu)=m$ and $\deg(\nu)=n$.
It is perhaps more convenient to consider $\Lambda$ as a graph
in which the vertices have degree $0$ and the edges are graded by
their (non-zero) degree, which takes values in $\bN^{\Bk}$, and
satisfy the unique factorization property above.

We are restricting our attention in this paper to $\Bk$-graphs with a single vertex.
Let $\ep_i$ for $1 \le i \le \Bk$ be the standard generators for $\bZ^\Bk$.
The generators of $\Fth$ are the maps of degree $\ep_i$, for $1 \le i \le \Bk$, 
which we label $e^i_\fs$ for $\fs \in \Bm_i = \{ 1,2,\dots,m_i\}$.
There are no commutation relations amongst the set $\{e^i_1,\dots , e^i_{m_i}\}$.
However the factorization property implies that each product $e^i_\fs e^j_\ft$
factors as $e^j_{\ft'}e^i_{\fs'}$.  The uniqueness of factorization implies that
there is a permutation $\theta_{ij}$ in $S_{\Bm_i \times \Bm_j}$ so that
\[
 e^i_\fs e^j_\ft = e^j_{\ft'} e^i_{\fs'} 
 \quad\text{where}\quad 
 \theta_{ij}(\fs,\ft) = (\fs',\ft') .
\]
The family $\{\theta_{ij}: 1 \le i < j \le \Bk \}$ is denoted by $\theta$;
and the $\Bk$-graph $\Fth$ is the semigroup generated by
$\{e^i_\ft : 1 \le i \le \Bk,\ 1 \le \ft \le m_i \}$.

Unfortunately, not every family of permutations $\theta$ yields a $\Bk$-graph.
There are evidently issues about associativity of the product and uniqueness
of the factorization.
For $\Bk=2$, every permutation yields a 2-graph;
but this is not true for $\Bk\ge3$.
See, for example, \cite{FowSim, RaeSimYee1}.
Fowler and Sims \cite{FowSim} showed that for $\Bk\ge 3$, 
$\theta$ determines a $\Bk$-graph if and only if 
every three sets of generators satisfy the following 
\textit{cubic condition} 
showing that a word of degree $(1,1,1)$ is well defined:
\begin{align*}
e^i_{\ft_1} e^j_{\ft_2} e^k_{\ft_3}
&= e^i_{\ft_1} e^k_{\ft_3'} e^j_{\ft_2'}
= e^k_{\ft_3''} e^i_{\ft_1'} e^j_{\ft_2'}
= e^k_{\ft_3''} e^j_{\ft_2''} e^i_{\ft_1''} \\
&= e^j_{\ft_{2'}} e^i_{\ft_{1'}} e^k_{\ft_3}
= e^j_{\ft_{2'}} e^k_{\ft_{3'}} e^i_{\ft_{1''}}
= e^k_{\ft_{3''}} e^j_{\ft_{2''}} e^i_{\ft_{1''}} 
\end{align*}
\[ \quad \text{implies}\quad 
e^i_{\ft_1''}=e^i_{\ft_{1''}},\
e^j_{\ft_2''}=e^j_{\ft_{2''}} \AND 
e^k_{\ft_3''}=e^k_{\ft_{3''}}.
\]
Hence, for $\Bk\ge 3$, $\Fth$ is a $\Bk$-graph if and only if the restriction of $\Fth$
to every triple family of edges $\{ e^i_\fs,\ e^j_\ft,\ e^k_\fu \}$
is a 3-graph.

We can consider each permutation $\theta_{ij}$ as a permutation
of $\prod_{i=1}^\Bk \Bm_i$ which fixes the coordinates except for
$i,j$, on which it acts as $\theta_{ij}$.  With this abuse of notation,
one can rephrase the cubic condition as:
\[
 \theta_{ij}\theta_{ik}\theta_{jk}=\theta_{jk}\theta_{ik}\theta_{ij}
 \qforal 1 \le i< j < k \le \Bk .
\]
This will facilitate calculations. 

We provide a few examples.

\begin{example}
Power \cite{Power1} showed that there are nine 2-graphs with $m_1=m_2=2$
up to isomorphism.  In \cite{DYperiod}, we showed that only two of these
are periodic (defined in the next section).  
These are the flip algebra in which $\theta(\fs,\ft)=(\ft,\fs)$ and
the square algebra given by the permutation
$\big( (1,1), (1,2), (2,2), (2,1) \big)$.

A more typical example of a 2-graph is the 
forward 3-cycle semigroup given by the
permutation $\big((1,1),(1,2),(2,1)\big)$. 
Curiously, the reverse $3$-cycle semigroup
arising from the $3$-cycle $\big((1,1),(2,1),(1,2)\big)$
is not isomorphic.
\end{example}

\begin{example}
\label{fff}
Let $m_i = n$ for all $1 \le i \le \Bk$
and $\theta_{ij}$ be the transposition
$\theta_{ij}(\fs,\ft)=(\ft,\fs)$.
Equivalently, this means that $e^i_\fs e^j_\ft = e^j_\fs e^i_\ft$
for all $i,j$ and all $\fs,\ft$.
It is readily calculated that 
\[
 \theta_{ij}\theta_{ik}\theta_{jk}(\fs,\ft,\fu) = (\fu,\ft,\fs) = 
 \theta_{jk}\theta_{ik}\theta_{ij}(\fs,\ft,\fu) .
\]
Thus this is a $\Bk$-graph.
\end{example}

\begin{example}
\label{ffflip}
Let $m_1=m_2=m_3=2$.
Let $\theta_{13} = \theta_{23}$ be the forward 3-cycle 
$\big((1,1),(1,2),(2,1)\big)$ and let $\theta_{12}$ be the flip.
Observe that $\theta_{13}(\fs,\ft) = (\ft,\fs+\ft)$ where addition is
calculated in $\bZ/2\bZ$.  Thus $\theta$ yields a 3-graph because
\[
 \theta_{12}\theta_{13}\theta_{23}(\fs,\ft,\fu) = (\fu,\ft+\fu,\fs+\ft+\fu) = 
 \theta_{23}\theta_{13}\theta_{12}(\fs,\ft,\fu) .
\]
\end{example}

\begin{example}
\label{ssf}
Let $m_1=m_2=m_3=2$.
Let $\theta_{13} = \theta_{23}$ be the square algebra 
which can be written $\theta_{13}(\fs,\ft)=(\ft,\fs+1)$, 
and let $\theta_{12}$ be the flip.
Then $\theta$ determines a 3-graph since
\[
 \theta_{12}\theta_{13}\theta_{23}(\fs,\ft,\fu) = (\fu,\ft+1,\fs+1) = 
 \theta_{23}\theta_{13}\theta_{12}(\fs,\ft,\fu) .
\]
\end{example}

\subsection{Representations}
Now consider the representations of $\Fth$.

A \textit{partially isometric representation} of $\Ftheta$ is a
semigroup homomorphism $\sigma: \Ftheta\to \B(\H)$ whose range
consists of partial isometries on a Hilbert space $\H$.
The representation $\sigma$ is  \textit{isometric} if the range
consists of isometries.

Call $\sigma$ \textit{row contractive} if  $[\sigma(e^i_1) \cdots \sigma(e^i_{m_i})]$
is a row contraction for $1 \le i \le \Bk$; and is
\textit{row isometric} if these row operators are isometries.
A row contractive representation is \textit{defect free} if
\[
 \sum_{\fs=1}^{m_i} \sigma(e^i_\fs) \sigma(e^i_\fs)^* = I
 \qforal 1 \le i \le \Bk .
\]
A row isometric defect free representation is a
\textit{$*$-representation} of $\Fth$.

The row isometric condition is equivalent to saying that the
$\sigma(e^i_\fs)$'s  are isometries with pairwise orthogonal range
for each $1 \le i \le \Bk$.

The most basic example of an isometric representation of $\Fth$
is the \textit{left regular representation $\lambda$}.  
This is defined on $\Fockth$ with orthonormal basis $\{ \xi_w : w \in \Fth\}$
given by $\lambda(v)\xi_w = \xi_{vw}$.
Each $\lambda(e^i_\fs)$ is an isometry.
Because the factorization of an element in $\Fth$ can begin with a unique
element of $\{e^i_1,\dots, e^i_{m_i}\}$ if it has any of these
elements as factors, it is clear that the ranges of $\lambda(e^i_\fs)$
are pairwise orthogonal for $1 \le \fs \le m_i$.
Hence this is a row isometric representation.
However it is also evident that it is not defect free.

One forms a $*$-algebra $\A$ generated by $\Fth$ subject to the 
relations implicit in $*$-representations that each $e^i_\fs$ is
an isometry and \mbox{$\sum_{\fs=1}^{m_i} e^i_\fs e^{i*}_\fs = 1$}
for $1 \le i \le \Bk$.  It is an easy exercise to see that $\A$
is the span of words of the form $uv^*$ for $u,v \in \Fth$.
This follows from the identity
\begin{align*}
 e^{i*}_\fs e^j_\ft &= e^{i*}_\fs e^j_\ft 
 \sum_r e^i_r e^{i*}_r  
 = \sum_r e^{i*}_\fs e^i_{r'} e^j_{\ft_r} e^{i*}_r
 = \sum_r \delta_{\fs r'} e^j_{\ft_r} e^{i*}_r.
\end{align*}
Every $*$-representation $\pi$ of $\Fth$ extends to a representation of $\A$.

The $\Bk$-graph C*-algebra $\ca(\Fth)$ is the universal C*-algebra
for $*$-representations of $\Fth$.
This is the completion of $\A$ with respect to the norm
\[ \|A\| = \sup \{ \| \pi(A)\| : \pi \text{ is a $*$-representation} \}.\]
This is the C*-algebra with the universal property that every
$*$-rep\-re\-sen\-tation of $\Fth$ extends uniquely to a $*$-representation
of $\ca(\Fth)$.

The universal property of $\ca(\Fth)$ yields a family of
\textit{gauge automorphisms}.
For any character $\phi$ in the dual group $\wh{\bZ^\Bk} \cong \bT^\Bk$ 
of $\bZ^\Bk$, consider the representation $\gamma_\phi(w) = \phi(\deg(w)) w$.
This is evidently an automorphism of $\Fth$.
So by the universal property of $\ca(\Fth)$, it extends to a $*$-automorphism
of $\ca(\Fth)$, which we also denote $\gamma_\phi$.

Integration over the $\bT^\Bk$ yields a faithful expectation
\[ \Phi(X) = \int_{\bT^\Bk}  \gamma_\phi(X) \,d\phi .\]
Checking this map on words $uv^*$, one readily sees that
$\Phi(uv^*) = \delta_{k0} uv^*$ where $k = \deg(uv^*) := \deg(u)-\deg(v)$.
Therefore 
\[ \fF:=\ran\Phi = \spn\{uv^* : \deg(uv^*)=0 \} .\]
Kumjian and Pask \cite{KumPask} show that this is an AF-algebra.
In our case of a single vertex, it is the UHF algebra for the
supernatural number $\prod_{i=1}^\Bk m_i^\infty$.
In particular, it is simple.

\subsection{Dilations}
If $\sigma$ is a representation of an operator algebra $\A$
on a Hilbert space $\H$, we say that a representation $\rho$ on a Hilbert
space $\K$ containing $\H$ is a dilation of $\sigma$ if $\sigma(a) = P_\H \rho(a)|_\H$
for all $a\in\A$.  This implies that $\H$ is semi-invariant, i.e. $\H = \M_1 \ominus\M_2$
for two invariant subspaces $\M_2 \subset \M_1$.
Arveson's dilation theory \cite{Arv1}, extended by Hamana \cite{Ham}, 
shows that $\A$ sits in a canonical C*-algebra known as its C*-envelope, $\cenv(\A)$.
This is determined by the universal property that whenever $j:\A \to \B(\H)$ is
a completely isometric isomorphism, there is a unique $*$-homomorphism
of $\ca(j(\A))$ onto $\cenv(\A)$ which extends $j^{-1}$.

A recent proof of Hamana's Theorem by Dritschel and McCullough \cite{DMc}
shows that every completely contractive representation of $\A$ dilates to a
maximal dilation, in the sense that any further dilation can only be obtained
by a direct sum of another representation.  Moreover, these maximal representations
of $\A$ are precisely those representations which extend to a $*$-representation
of $\cenv(\A)$.

The usefulness to us of dilations lies in the fact that a $*$-dilation
may be uniquely determined by a smaller, more tractable representation.
Indeed, we will be looking for a compression to a cyclic, coinvariant subspace
with some special structure.

The operator algebra that figures here is the nonself-adjoint unital operator 
algebra $\Ath$ generated by $\lambda(\Fth)$ as a subalgebra of $\B(\Fockth)$.
There is no simple criterion on a representation of $\Fth$ which is equivalent
to it having a completely contractive extension to $\Ath$.
A necessary condition is that the representation be row contractive because
\[
 \big\| \big[ \sigma(e^i_1)\ \dots \ \sigma(e^i_{m_i}) \big] \big\|
 \le \big\| \big[ \lambda(e^i_1)\ \dots \ \lambda(e^i_{m_i}) \big] \big\| = 1.
\]
However it was shown in \cite{DPYdilation} that this is a strictly weaker
condition than being completely contractive even for 2-graphs on one vertex.

Two results provide the information that we need, and they are both valid for
arbitrary $\Bk$-graphs, not just the single vertex case.
The first is a result of Katsoulis and Kribs  \cite{KK} on the C*-envelope of
higher rank graph algebras.

\begin{thm}[Katsoulis--Kribs] \label{T:cenv}
The C*-envelope of $\Ath$ is $\ca(\Fth)$.
\end{thm}

This implies that every completely contractive representation of $\Ath$
dilates to a $*$-representation of $\ca(\Fth)$, and hence to a $*$-representation of $\Fth$.

In \cite{DPYdilation}, we established a dilation theorem for a class of
doubly generated operator algebras which includes the 2-graphs on one vertex.
We showed, in particular, that defect free representations dilate to $*$-representations;
and consequently, they yields completely contractive representations of $\Ath$.
Using the Poisson transform defined by Popescu in \cite{Popescu},
Skalski and Zacharias \cite{SZ} studied the dilation theory of higher rank graphs
in a very general context. In particular, their results include the following dilation theorem
which is valid for all $\Bk$.

\begin{thm}[Skalski-Zacharias] \label{T:dilation}
Every defect free, row contractive representation of $\Fth$
has a unique minimal $*$-dilation.
\end{thm}

Consequently, every defect free representation of $\Fth$ extends to
a completely contractive representation of $\Ath$.

The algebra $\Ath$ will not play a direct role in the current paper, which is focussed
on $*$-representations and the structure of $\ca(\Fth)$.  However, it lurks in the
background because we use dilation theory to simplify the analysis of the
representations.

\section{Atomic Representations}

Atomic representations of 2-graphs are comprehensively studied in \cite{DPY}.
They provide a very interesting class of representations,
and play an important role in the study of 2-graphs. 
We now introduce such representations of $\Bk$-graphs.

\begin{defn}
A partially isometric representation is \textit{atomic}
if there is an orthonormal basis which is permuted,
up to scalars, by each partial isometry.
That is, $\sigma $ is atomic if there is an orthonormal  basis  \mbox{$\{\xi_k : k\ge1\}$}
so that for each $w \in \Fth$, $\sigma(w) \xi_k = \alpha \xi_l$ for
some $l$ and some $\alpha \in \bT \cup \{0\}$.
\end{defn}

As in \cite{DPY},
we refer to the atomic partially isometric representations in this paper
simply as atomic representations, and
likewise we refer to the row contractive defect free
atomic representations simply as defect free atomic representations.

Every atomic representation determines a graph with the standard basis vectors
representing the vertices, and the partial isometries $\sigma(e^i_\fs)$ determining
directed edges: when $\sigma(e^i_\fs) \xi_a = \xi_b$, we draw an edge from vertex $a$
to vertex $b$ labelled $e^i_\fs$.  The graph is called connected if there is an undirected
path from each vertex to every other. This graph contains all information about
the representation except for the scalars of modulus one. 

As in \cite[Lemma~5.6]{DPYdilation}, it is easy to see that the $*$-dilation of an atomic
defect free representation remains atomic. We leave the straightforward extension
to $\Bk$-graphs to the reader.

\begin{lem} \label{L:dilate atomic}
If $\sigma$ is an atomic defect free representation, 
with respect to some basis of $\H_\sigma$, with minimal 
$*$-dilation $\pi$, then this basis extends to a basis
for $\H_\pi$ making $\pi$ an atomic $*$-representation.
\end{lem}

\subsection{Inductive Limit Representations}

In this subsection, we will define a whole family of $*$-dilations of $\lambda$
which are, in fact, inductive limits of $\lambda$.
They play a central role in what follows.

Arbitrarily choose an \textit{infinite tail} $\tau$ of $\Fth$;
that is, an infinite word in the generators which has infinitely
many terms from each family $\{e^i_1,\dots, e^i_{m_i}\}$.
Such an infinite word can be factored so that these terms occur
in succession as
\[
\tau=e^1_{\ft_{01}}\cdots e^\Bk_{\ft_{0\Bk}}
     e^1_{\ft_{11}}\cdots e^\Bk_{\ft_{1\Bk}}\cdots 
     = \tau_0\tau_1\tau_2\dots
\]
where $\tau_s = e^1_{\ft_{s1}}\cdots e^\Bk_{\ft_{s\Bk}}$ for $s \ge0$.
Let $\F_s = \F := \Ftheta$ for $s\ge0$, viewed as discrete
sets on which the generators of $\Ftheta$ act as injective maps
by right multiplication, namely, $\rho(w)f = fw$ for all  $f \in \F$.
Consider $\rho_s = \rho(\tau_s)$ 
as a map from $\F_s$ into $\F_{s+1}$.
Define $\F_\tau$ to be the injective limit set
\[
 \F_\tau = \lim_{\rightarrow} (\F_s, \rho_s ) .
\]
Also let $\iota_s$ denote the injections of $\F_s$ into $\F_\tau$.
Thus $\F_\tau$ may be viewed as the union of $\F_0, \F_1, \dots $
with respect to these inclusions.

The left regular action $\lambda$ of $\Fth$ on itself induces
corresponding maps on $\F_s$ by  $\lambda_s(w) f = wf$.
Observe that $\rho_s \lambda_s = \lambda_{s+1} \rho_s$.
The injective limit of these actions is an action $\lambda_\tau$
of $\Fth$ on $\F_\tau$.
Let $\lambda_\tau$ also denote the corresponding representation of $\Fth$
on $\ltwo(\F_\tau)$.
That is, we let  $\{ \xi_f : f \in \F_\tau\}$ denote the orthonormal basis
and set $\lambda_\tau(w) \xi_f = \xi_{wf}$.
It is easy to see that this provides a defect free,
isometric representation of $\Fth$; i.e.\ it is a $*$-representation.

Let $\tau^s = \tau_0\tau_1\dots \tau_s$ for $s \ge 0$.
We may consider an element $w = \iota_s(v)$
as $w = v\tau^{s*}$.
This makes sense in that 
$\xi_w = \lambda_\tau(v) \lambda_\tau(\tau^s)^* \xi_{\iota_0(\mt)}$.
In particular, we have $\xi_\mt =  \xi_{\iota_0(\mt)}$.
\medbreak

Since $\lambda_\tau$ is a $*$-dilation, it extends to a representation 
of $\ca(\Fth)$ which we also denote by $\lambda_\tau$.
This is always a faithful representation.  
This is the analogue of \cite[Theorem~3.6]{DPYdilation}.

\begin{thm} \label{T:faithful}
For any infinite tail $\tau$, the representation $\lambda_\tau$
is faithful on $\ca(\Fth)$.
\end{thm}

\begin{proof}
Because of the gauge invariance uniqueness theorem, it suffices
to show that there are gauge automorphisms of $\ca(\lambda_\tau(\Fth))$.
This is accomplished by conjugation by a diagonal unitary.
Given a word $w = v\tau^{s*} \in \F_\tau$, define $\deg(w) = \deg(v) - \deg(\tau^s)$.
It is clear that this is well defined and extends the degree map on $\Fth$.
Given $\phi \in \wh{\bZ^\Bk}$,  define 
$U_\phi = \diag\big( \phi(\deg(w)) \big)$ 
with respect to the basis $\{\xi_w : w \in \F_\tau\}$.
Then for $u \in \Fth$ and $w= v\tau^{s*} \in \F_\tau$,
\begin{align*}
 U_\phi \lambda_\tau(u) U^*_\phi \xi_w &=
 U_\phi \lambda_\tau(u) \ol{\phi(\deg(w))} \xi_w \\
 &= U_\phi \ol{\phi(\deg(v)-\deg(\tau^s))} \xi_{uw}\\
 &= \phi(\deg(uv)-\deg(\tau^s)) \ol{\phi(\deg(v)-\deg(\tau^s))} \xi_{uw} \\
 &= \phi(\deg(u)) \lambda_\tau(u) \xi_w .
\end{align*}
Hence $U_\phi \lambda_\tau(u) U^*_\phi = \lambda_\tau(\gamma_\phi(u))$.
\end{proof}

Because of the dilation theory, $\lambda_\tau$ is completely
determined by any cyclic coinvariant subspace $\H$ of $\ltwo(\F_\tau)$
as the unique minimal $*$-dilation of this compression.
We will describe such a subspace which will be convenient.

Observe that $\lambda_\tau(w)^* \xi_{\iota_0(\mt)}$ is non-zero if and
only if $w$ is an initial segment of $\tau$ after appropriate factorization.
That is, given $n = (n_1,\dots,n_\Bk) \in \bN_0^\Bk$, one may factor
$\tau$ in exactly one way so that $\tau =w_n \tau_n'$ in which
$w_n$ has degree $n$.
Let $\zeta_{-n} := \lambda_\tau(w_n)^* \xi_{\iota_0(\mt)}$ 
for $n \in \bN_0^\Bk$.  In particular, $\zeta_{(-s,\dots,-s)} = \xi_{\iota_s(\mt)}$.
Then 
\[
 \H_\tau = \spn\{ \zeta_n : n \in (-\bN_0)^\Bk \}
\]
is evidently a cyclic subspace because it contains each $\xi_{\iota_s(\mt)}$,
and is coinvariant by construction.

Note that beginning at any of the standard basis vectors $\xi_w$, there will be
some word $v$ so that $\lambda_\tau(v)^* \xi_w$ is a basis vector $\zeta_n$ in $\H_\tau$.
As the restriction of $\lambda_\tau$ to the cyclic subspace generated by $\zeta_n$ is
unitarily equivalent to $\lambda$, it is easy to understand why the restriction
of $\lambda_\tau$ to $\H_\tau$ determines the whole representation.

For each $n \in -\bN_0^\Bk$, there are unique integers $\ft_n^i$
so that $\zeta_n$ is in the range of $\lambda_\tau(e^i_{\ft_n^i})$
for $1 \le i \le \Bk$; that is, $\lambda_\tau(e^i_{\ft_n^i}) \zeta_{n-\ep_i} = \zeta_n$. 
Set $\Sigma(\tau,n) = (\ft^1_{n},...,\ft^\Bk_{n}) $.
This determines the data set
\[
 \Sigma(\tau) = 
 \{ \Sigma(\tau,n) : n \in -\bN_0^\Bk \} .
\]

\begin{defn}
Two tails $\tau_1$ and $\tau_2$ 
are said to be \textit{tail equivalent} if their data sets eventually coincide;
i.e.\ there is $T \in -\bN_0^\Bk$ so that
\[
  \Sigma(\tau_1,n) =  \Sigma(\tau_2,n)  \qforal n \le T .
\]
We say that $\tau_1$ and $\tau_2$ are \textit{$p$-shift tail equivalent}
for some $p \in \bZ^k$ if there is a $T\in\bN_0^k$ so that
\[
  \Sigma(\tau_1,n) =  \Sigma(\tau_2,n+p)  \qforal n \le T .
\]
Then $\tau_1$ and $\tau_2$ are \textit{shift tail equivalent} if they
are $p$-shift tail equivalent for some $p\in \bZ^\Bk$.
\end{defn}

Clearly, if $\tau_1$ and $\tau_2$ are shift tail equivalent, then $\lambda_{\tau_1}$
and $\lambda_{\tau_2}$ are unitarily equivalent.

We now introduce two important notions: the symmetry of $\tau$
and the aperiodicity condition of $\Fth$.

\begin{defn}
A tail $\tau$ is \textit{$p$-periodic} if 
$ \Sigma(\tau,n) =  \Sigma(\tau,n+p)$ for all $n \le 0 \wedge -p$;
and  \textit{eventually $p$-periodic} if $\tau$  is $p$-shift tail equivalent to itself.
The \textit{symmetry group} of $\tau$  is the subgroup of $\bZ^\Bk$ 
\[
 H_\tau = \{ p \in \bZ^\Bk :
 \tau \text{ is eventually $p$-periodic} \} .
\]
The \textit{symmetry group} of $\Fth$ is defined by
\[
H_\theta :=\cap_\tau H_\tau
\]
as $\tau$ runs over all possible infinite tails of $\Fth$.

A tail $\tau $ is called \textit{aperiodic} if $H_\tau = \{0\}$.
The semigroup $\Fth$ is \textit{aperiodic} if $H_\theta=\{0\}$.
Otherwise we say that $\Fth$ is \textit{periodic}.
\end{defn}

Clearly if there is an aperiodic infinite tail $\tau$, then $\Fth$ is aperiodic.
The following result shows that our definition coincides with the 
Kumjian--Pask aperiodicity condition.

\begin{prop} \label{P:goodtail}
$\Fth$ has an infinite tail with $H_\tau = H_\theta$.
In particular, when $\Fth$ is aperiodic, there is an aperiodic tail.
More generally, we have \mbox{$H_\theta \cap \bN_0^\Bk = \{0\}$.}
\end{prop}

\begin{proof}
For each $p$ in $\bZ^\Bk \bsl H_\theta$, there is a tail $\tau$
such that $p \not\in H_\tau$.  Hence there is some $n \in -\bN_0^\Bk$
so that $n+p \in -\bN_0^\Bk$ and $\Sigma(\tau,n) \ne \Sigma(\tau,n+p)$.
Choose $s$ so that $0 \ge n, n+p \ge (-s,\dots,-s)$.
Then the finite initial segment $w_p=\tau_0\dots\tau_s$ of $\tau$ 
already exhibits the lack of $p$-symmetry.
So any infinite tail that contains $w_p$ infinitely often can never
exhibit $p$-symmetry for $n \le T$.
Form a tail $\tau$ by stringing together the words $w_p$, repeating each
one infinitely often.  By construction, $H_\tau \subset H_\theta$.
The other inclusion is true by definition.

If $p \in \bN_0^\Bk \bsl\{0\}$, it is easy to write down a finite sequence
which does not have $p$-symmetry.  Splicing such words into $\tau$
as above shows that $H_\theta \cap \bN_0^\Bk = \{0\}$.
\end{proof}

\section{A group construction}
\label{S:Groupcons}

In this section, we describe a large family of defect free atomic representations
with a very nice structure.  By the dilation theorem, they encode a
family of atomic $*$-representations which are obtained as the unique minimal
$*$-dilations.  The main result is that every irreducible atomic $*$-representation
arises from this construction; and every atomic
$*$-representation decomposes as a direct integral of the irreducible ones.

Let $G$ be a finitely generated abelian group with $\Bk$ designated 
generators $\fg_1,...,\fg_\Bk$.
Suppose that functions are given 
\begin{alignat*}{2}
 \ft^i &:  G\to \{1,...,m_i\}, \quad
                &\ft^i(g)    &=: \ft^i_g ,\quad i=1,...,\Bk, \\
 \alpha^i &:  G\to \bT,  \quad
                &\alpha^i(g) &=: \alpha^i_g ,\quad i=1,...,\Bk.
\end{alignat*}

Consider a defect free atomic representation $\sigma:\Fth\to\B(\ell^2(G))$ 
given by
\[
 \sigma(e^i_\ft)\xi_{g-\fg_i} =\delta_{\ft, \ft^i_g}\, \alpha^i_g\,  \xi_g
 \qfor i=1,...,\Bk.
\]
This ensures that $\xi_g$ is in the range of 
$\sigma(e^i_{\ft^i_g})$ for each $1 \le i \le \Bk$.
In order for $\sigma$ to be a representation, the commutation relations
must be satisfied, namely
\[
 e^i_{\ft_g} e^j_{\ft_{g-\fg_i}} =
 e^j_{\ft_g} e^i_{\ft_{g-\fg_j}} 
 \qand
 \alpha^i_{\ft_g} \alpha^j_{\ft_{g-\fg_i}} = 
 \alpha^j_{\ft_g} \alpha^i_{\ft_{g-\fg_j}} 
\]
for all $g \in G$ and $1 \le i < j \le \Bk$. 
Such a representation will be called a \textit{group construction representation}.

\begin{eg}
Actually defining such relations might not be so easy.
However in the case of $G = \bZ^\Bk$ with the standard generators,
we can obtain all such representations from infinite tails.
Indeed, we saw that an infinite tail $\tau$ gives rise to an
inductive limit representation $\lambda_\tau$.
We then identified a subspace $\H_\tau$ with basis $\{\zeta_n : n \in -\bN_0^\Bk\}$.
It is possible to continue this `forward' to obtain a (non-canonical)
representation modelled on the group $\bZ^\Bk$.

This may be accomplished by extending $\tau$ arbitrarily to a doubly infinite word,
say $\omega^*\tau$, where $\omega^*$ is a tail in reverse order.
This will specify how we are allowed to move forward and stay
within our subspace.
If $n \in -\bN^\Bk$, recall that there is a unique word $w_n$ so that
$\tau$ factors as $\tau = w_n \tau'$.  Suppose that $m \in \bN_0^\Bk$.
Then there is a unique word $v_m$ of degree $m$ so that 
$\omega^* = \omega^{\prime*} v_m$.
Let $\zeta_{n+m} = \lambda_\tau(v_mw_n) \zeta_n$.
It is not difficult to verify that this is well defined.

We identify $\K = \spn\{\zeta_m : m \in \bZ^\Bk\}$ with $\ltwo(\bZ^\Bk)$.
It is easy to see that $\K$ is a coinvariant subspace; and it is cyclic
because $\H_\tau$ is cyclic.
Let $\sigma$ be the compression of $\lambda_\tau$ to this subspace.
Then we have a representation of group type with all constants $\alpha^i_g = 1$.
\end{eg}

It turns out that the scalars $\alpha^i_g$ are not difficult to control. 
It was shown in \cite[Theorem~6.1]{DPY} that the
representation is unitarily equivalent to another group construction 
representation in which these constants are independent of $g \in G$.
The proof is essentially the same as the 2-graph case.  
So we state it without proof.

For each group $G$, there is a canonical homomorphism $\kappa$
of $\bZ^\Bk$ onto $G$ sending the standard generators 
$\ep_i$ to $\fg_i$ for $1 \le i \le \Bk$.
Let $K$ be the kernel of $\kappa$, so that $G \cong \bZ^\Bk/K$.

\begin{thm} \label{T:scalars}
Let $G = \bZ^\Bk/K$ as above, and let $\sigma$ be a group construction
representation.  Then $\sigma$ is unitarily equivalent to another group construction
representation with the same functions $\ft^i$ and constant functions $\alpha^i$.
\end{thm}

In fact, as for 2-graphs, the constants determine a unique character 
$\psi$ of $K$.  If a path in the graph is a loop returning to the vertex
where it started, then it determines a unique word in the generators and
their adjoints whose degree $d$ belongs to $K$.  The scalar multiple of the
vertex vector obtained by application of the partial isometry is $\psi(d)$.
The choice of the scalars $\alpha^i$ is determined by choosing
an extension $\phi$ of $\psi$  to a character on $\bZ^\Bk$; and they are 
given by $\alpha_i = \phi(\ep_i)$.
The character $\phi$ is determined by $\psi$ and a character $\chi$ of $G$;
and the constants $\alpha^i$ can be replaced by $\alpha^i \chi(\fg_i)$.
We will not need this detailed information, so we refer the interested reader
to the proof in \cite{DPY} in the 2-graph case.

There are two useful notions of symmetry for these group constructions.

\begin{defn}
A group construction representation $\sigma$ of $\Fth$ on $\ltwo(G)$
has a \textit{full symmetry subgroup} $H \le G$ if $H$ is the largest
subgroup of $G$ such that
$\ft^i_{g+h}=\ft^i_g$ and $\alpha^i_{g+h} = \alpha^i_g$
for all $g \in G$ and $h \in H$ and $1 \le i \le \Bk$. 

We define the
\textit{symmetry group} $H_\sigma \le G$ to be the largest subgroup
$H \le G$ for which there is some $T=(T_1,\dots,T_\Bk) \in \bZ^\Bk$ such that
$\ft^i_{g+h}=\ft^i_g$ and $\alpha^i_{g+h} = \alpha^i_g$
for all $g \in G$ and $h \in H$ and $1 \le i \le \Bk$ with both 
$g,g+h \in G_T = \big\{ \sum n_i \fg_i : n_i \le T_i \big\}$.

We also will say that two group constructions $\sigma$ and $\tau$
on $\ltwo(G)$ given by functions $\fs^i_g$, $\alpha^i_g$ and $\ft^i_g$, $\beta^i_g$
respectively are \textit{tail equivalent} if 
$\ft^i_g =\fs^i_g$ and $\beta^i_g = \alpha^i_g$
for all $g \in G_T$ and $1 \le i \le \Bk$.
\end{defn}

Note that the symmetry group $H_\sigma$ is well defined because an
increasing sequence of subgroups of a finitely generated abelian group $G$
is eventually constant.

The significance of tail equivalence is that two tail equivalent representations
will have unitarily equivalent $*$-dilations.  Indeed, by the Skalski--Zacharias
dilation theorem, each has a unique $*$-dilation.  However, both must
coincide with the unique $*$-dilation of their common restriction to $\ltwo(G_T)$.

We wish to show that these group constructions may be obtained
by dilating certain partial group constructions, preserving the
symmetry.  In particular, we will show that a group construction $\sigma$
with symmetry group $H$ is tail equivalent to one with full symmetry group $H$.

\begin{thm} \label{T:dilate to group}
Let $G$ be an abelian group with generators $\fg_1,\dots,\fg_\Bk$.
Given $T=(T_1,\dots,T_\Bk)$ with $T_i \in \bZ \cup \{\infty\}$,
let \mbox{$G_T = \big\{\! \sum n_i \fg_i : n_i \le T_i \big\}$.}
Suppose that $\sigma$ is a representation of $\Fth$ on $\ltwo(G_T)$,
determined by functions $\ft^i:G_T \to \Bm_i$ and $\alpha^i : G_T \to \bT$ so that
\[
 \sigma(e^i_\ft)\xi_{g-\fg_i} =\delta_{\ft, \ft^i_g}\, \alpha^i_g\,  \xi_g
 \qfor i=1,...,\Bk \AND g \in G_T.
\]
Then $\sigma$ may be dilated to a group construction representation on $\ltwo(G)$.
\end{thm}

\begin{proof}
Consider dilations to $\ltwo(G_{T'})$ for $T' \ge T$.
Among these dilations, select one which is maximal in the sense that it cannot
be dilated to a representation on a larger subset of this type.
There is no loss of generality in assuming that this set is $G_T$ itself.
We will show that $G_{T}=G$.

Indeed, otherwise there is some $T_i < \infty$.
By the Skalski--Zacharias dilation theorem, there is a $*$-dilation $\tilde\sigma$ of $\sigma$
on a Hilbert space $\H$ containing $\ltwo(G_T)$, and by Lemma~\ref{L:dilate atomic}, this
representation is atomic.

Assume first that every $T_i < \infty$.  Set $a = \sum_{i=1}^\Bk T_i \fg_i$; 
and arbitrarily pick some $i$.  Let $T'=T+\ep_i$, so
\[
 G_{T'} = G_T + \fg_i =  G_T \cup (a+\fg_i+G_0) = G_T \dcup (a+\fg_i + S)
\]
where  $S = \big\{ \sum_{j \ne i} n_j \fg_j : n_j \le 0 \}$.
We will extend $\sigma$ to $\ltwo(G_{T'})$.
To this end, arbitrarily select $p_a \in \Bm_i$.
Identify $\xi_{a+\fg_i}$ with $\tilde\sigma(e^i_{p_a})\xi_a$ in $\H$;
and set $\ft^i_a = p_a$ and $\alpha^i_a = 1$.
For each $s = \sum_{j \ne i} n_j \fg_j$ in $S$, $a+s \in G_T$ and 
there is a unique word $w \in \Fth$ of degree $(|n_1|,\dots,|n_\Bk|)$
so that $\sigma(w) \xi_{a+s} = \alpha \xi_a$ for some $\alpha \in \bT$.
Factor $e^i_{p_a}w = w' e^i_{p_{a+s}}$.  
Then identify $\xi_{a+s+\fg_i}$ with $\tilde\sigma(e^i_{p_{a+s}})\xi_{a+s}$ in $\H$;
and set $\ft^i_{a+s} =p_{a+s}$ and $\alpha^i_{a+s}=1$.

When $n_j \ne 0$, we need to define $\ft^j_{a+s+\fg_i}$.
Observe that $a+s+\fg_j \in G_T$.  
If $\ft^j_{a+s} = q$, we can factor $w = ve^j_q$.
Now 
\[ \sigma(w) \xi_{a+s} = \sigma(v) \alpha^j_{a+s} \xi_{a+s+\fg_j} = \alpha \xi_a .\]
Factor $e^i_{p_a}v = v' e^i_{p_{a+s+\fg_j}}$ as before;
and factor $e^i_{p_{a+s+\fg_j}} e^j_q = e^j_{q'} e^i_{p'}$.
Then we have 
\[
 v' e^j_{q'} e^i_{p'} = v' e^i_{p_{a+s+\fg_j}} e^j_q = 
 e^i_{p_a} v e^j_q = e^i_{p_a}w = w' e^i_{p_{a+s}} = 
 v' e^j_{q'} e^i_{p_{a+s}} .
\]
Therefore $p'=p_{a+s}$ and
\begin{align*}
 \tilde\sigma( e^j_{q'} ) \xi_{a+s+\fg_i} &= 
 \tilde\sigma(e^j_{q'}e^i_{p_{a+s}})\xi_{a+s} =
 \tilde\sigma(e^i_{p_{a+s+\fg_j}} e^j_q)\xi_{a+s} \\&= 
 \tilde\sigma(e^i_{p_{a+s+\fg_j}}) \alpha^j_{a+s}\xi_{a+s+\fg_j} =
 \alpha^j_{a+s} \xi_{a+s+\fg_j+\fg_i}
\end{align*}
So we set $\ft^j_{a+s+\fg_i} = q'$ and $\alpha^j_{a+s+\fg_i} = \alpha^j_{a+s}$.
It now follows that we have defined a dilation of $\sigma$ to $\ltwo(G_{T'})$.
This contradicts the maximality of $G_T$.

Now consider the case in which $T_j=\infty$ for $j \in J$, 
for some subset $J \subset \{1,\dots,\Bk\}$.  
If $G_T \ne G$, then there is some $i$ with $T_i < \infty$.
Let $a_k \in \bN_0^\Bk$ where $a_k^i = T_i$ when $T_i<\infty$ and $a_k^j = k$ when $j \in J$.
Arguing as above, we can construct a dilation $\sigma_k$ to $\ltwo(A_k)$
where $A_k = G_T \dcup (a+\fg_i+S)$.
If we look at the action of $\sigma_k$ at $\xi_{a_l}$ for $0 \le l \le k$, then
one sees that the value of $p_{a_l} = \ft^i_{a_l}$ takes some value infinitely often.
Using a diagonal argument, we may drop to a subsequence so that
the values of $p_{a_l}$ are constant for a sequence $\sigma_{k_s}$.
Hence it is apparent that one can define a representation $\sigma'$
on $\ltwo(A_\infty)$, where $A_\infty = \bigcup_k A_k = G_T + \fg_i = G_{T'}$
which extends $\sigma$.
Since $G_T$ was presumed to be maximal, we obtain $G_T=G$ as desired.
\end{proof}

\begin{cor} \label{C:dilate with symmetry}
Let $G$ be an abelian group with generators $\fg_1,\dots,\fg_\Bk$.
Suppose that $\sigma$ is a representation of $\Fth$ on $\ltwo(G_T)$
determined by functions $\ft^i:G_T \to \Bm_i$ and $\alpha^i : G_T \to \bT$ satisfying
\[
 \sigma(e^i_\ft)\xi_{g-\fg_i} =\delta_{\ft, \ft^i_g}\, \alpha^i_g\,  \xi_g
 \qfor i=1,...,\Bk \AND g \in G_T.
\]
Moreover, suppose that $H$ is a subgroup of $G$ such that
\[
 \ft^i_g = \ft^i_{g+h} \AND \alpha^i_g=\alpha^i_{g+h} 
 \quad\text{whenever}\quad g, g+h \in G_T.
\]
Then $\sigma$ may be dilated to a group construction representation on $\ltwo(G)$
with full symmetry group $H$.
\end{cor}

\begin{proof}
Using the symmetry group $H$, it is routine to extend the representation
to $\ltwo(G_T+H)$ with the same symmetry.  Then by collapsing the cosets of $H$
to single points, we obtain a representation $\tau$ of $G/H$ on $\ltwo((G_T+H)/H)$.
By Theorem~\ref{T:dilate to group}, there is a dilation $\tau'$ of $\tau$ on $\ltwo(G/H)$.
Unfolding this yields a representation on $\ltwo(G)$ with full symmetry group $H$
which dilates $\sigma$.
\end{proof}

\begin{cor} \label{C:symmetry group}
Suppose that $\sigma$ is a group construction representation of $\Fth$ on $\ltwo(G)$
with symmetry group $H$.  Then $\sigma$ is tail equivalent to another group
construction representation $\tau$ which has full symmetry group $H$.
\end{cor}

\begin{proof}
There is a $G_T$ so that $\sigma$ compressed to $\ltwo(G_T)$ has symmetry group $H$.
The previous corollary shows how to dilate this to $\tau$ with full symmetry group $H$.
Evidently, $\tau$ is tail equivalent to $\sigma$.
\end{proof}

\begin{cor} \label{C:inductive symmetry}
Suppose that $\tau$ is an infinite tail with symmetry group $H < \bZ^\Bk$.
Then there is a group construction representation of $\Fth$ on $\ltwo(\bZ^\Bk)$
with symmetry group $H$ with minimal $*$-dilation $\lambda_\tau$.
\end{cor}

\begin{proof}
By definition of the symmetry group of $\tau$, there is a coinvariant subspace
$\H$ spanned by $\{\zeta_n : n \le T\}$ which has $H$ symmetry.
By Corollary~\ref{C:dilate with symmetry}, this dilates to a group construction
$\sigma$ on $\bZ^\Bk$ with full symmetry group $H$.
The minimal $*$-dilation of $\sigma$ is a minimal $*$-dilation of 
$P_\H\lambda_\tau|_{\H}$---which is unique; and hence this dilation is $\lambda_\tau$.
\end{proof}

\section{Decomposing Atomic Representations}\label{S:Decomp}

We are now prepared to prove that every atomic $*$-representation may be decomposed
as a direct sum of $*$-representations obtained from dilating the group construction.
In \cite{DPY}, this was established for rank 2 graphs by a detailed case by case analysis.
So our new proof provides insight into that case as well.

 It is evident that the span of basis vectors
in any connected component of the graph is a reducing subspace; and every 
atomic representation decomposes as a direct sum of connected ones.  So in our
analysis of atomic representations, it suffices to consider the case of a connected graph.

\begin{prop} \label{P:unique tail}
Let $\sigma$ be an atomic $*$-representation with connected graph.
For any standard basis vector $\eta$, there is a unique infinite tail
$\tau = \tau_0\tau_1\tau_2\dots$ such that $\eta$ is in the range of
$\sigma(\tau_0\dots\tau_n)$ for all $n\ge0$.
The shift-tail equivalence class of $\tau$ is independent of the choice of $\eta$.
\end{prop}

\begin{proof}
Since each basis vector is in the range of $\sigma(e^i_j)$ for exactly one $j \in \Bm_i$,
the existence and uniqueness of $\tau$ is clear.
For any two basis vectors $\eta_1$ and $\eta_2$, there is a basis vector $\eta$
and words $u_1$ and $u_2$ so that $\sigma(u_i)\eta = \eta_i$.
If $\tau$ is the infinite tail obtained for $\eta$, it follows that the tails $\tau_i$ obtained 
for $\eta_i$ are $\tau_i = u_i \tau$.  So they are $\deg(u_2)-\deg(u_1)$ shift tail equivalent.
\end{proof}

As a consequence, we can define the symmetry group of $\sigma$ in an unambiguous way.

\begin{defn}
Define the \textit{symmetry group} of an atomic $*$-rep\-re\-sent\-ation with connected graph
as the symmetry group $H_\tau$ of any tail $\tau$ derived from any standard basis vector.
\end{defn}

\begin{thm} \label{T:atomic from group}
Every atomic $*$-representation $\sigma$ of $\Fth$ with a connected graph is obtained 
as the minimal $*$-dilation of a group construction representation $\rho$.
Moreover, the symmetry group of $\sigma$ coincides with 
the full symmetry group of $\rho$.
\end{thm}

\begin{proof}
Let $\sigma$ be an atomic $*$-representation with connected graph and standard basis
$\{\eta_j : j \ge0 \}$.
Start with an arbitrary basis vector, say $\eta_0$.
There is a unique infinite tail $\tau = \tau_0\tau_1\tau_2 \dots$ with
the degree of each $\tau_i$ equal to $(1,1,\dots,1)$ such that
$\eta_0$ is in the range of $\sigma(\tau_0\dots\tau_n)$ for all $n \ge 0$.
Let $\lambda_\tau$ be the inductive limit representation determined by $\tau$.
The standard basis for $\ltwo(\F_\tau)$ will be denoted $\{\xi_w : w \in \F_\tau \}$.

Let the symmetry group of $\lambda_\tau$ be $H < \bZ^\Bk$.
By Corollary~\ref{C:inductive symmetry}, there is a group construction
representation $\mu$ on $\ltwo(\bZ^\Bk)$ with full symmetry group $H$
which has $\lambda_\tau$ as its unique minimal $*$-dilation.
The basis $\{\zeta_n : n \in \bZ^\Bk\}$ is identified with a subset of the
basis of $\ltwo(\F_\tau)$ generating a coinvariant subspace.

There is a canonical map $\theta$ of $\{ \xi_w : w \in \F_\tau\}$ 
onto the basis $\{\eta_j  \}$
which intertwines the actions of $\lambda_\tau$ and $\sigma$ in the sense
\[ \sigma(v) \theta(\xi_w) = \theta(\xi_{vw}) \qforal v\in\Fth \AND w \in \F_\tau .\]
Indeed, every element of $\F_\tau$ may be factored as $w = u\tau^{s*}$
for some $s \ge0$ and $u \in \Fth$.  
Then we define $\theta(\xi_w) = \sigma(u)\sigma(\tau^s)^* \eta_0$.
In particular, $\theta$ carries $\{\zeta_n : n \in \bZ^\Bk \}$ 
onto the basis of a coinvariant subspace of $\sigma$.

Suppose that $\eta_j = \theta(\zeta_m) = \theta(\zeta_n)$ for $n,m \in \bZ^\Bk$.
Set $l =m \wedge n \in \bZ^\Bk$.  Then there are unique words $u,v \in \Fth$
of degrees $m-l$ and $n-l$ respectively so that 
\[
 \lambda_\tau(u) \zeta_l = \zeta_m \qand  \lambda_\tau(v) \zeta_l = \zeta_n .
\]
Hence $\lambda_\tau(u)\lambda_\tau(v)^* \zeta_n = \zeta_m$.
Therefore $\sigma(u)\sigma(v)^* \eta_j = \eta_j$.
Conversely, if $\sigma(u)\sigma(v)^* \eta_j = \eta_j$ and $\theta(\zeta_n) = \eta_j$,
then $\theta(\zeta_m) = \eta_j$ where $m-n = \deg(u) - \deg(v)$.

Observe that $K_j = \{m-n : \eta_j = \theta(\zeta_m) = \theta(\zeta_n) \}$ 
is a subgroup of $H$. Indeed, it is clear that if $\theta(\zeta_m)=\theta(\zeta_n)$,
then the infinite tails obtained by pulling back from the vectors $\zeta_m$
and $\zeta_n$ both coincide with the tail obtained by pulling back from $\eta_j$.
Hence $m-n \in H$.  Consequently, if $\deg(uv^*) = m-n$ such that 
$\lambda_\tau(u) \lambda_\tau(v)^* \zeta_n = \zeta_m$ and $\theta(\zeta_l) = \eta_j$, 
it follows that $\lambda_\tau(u) \lambda_\tau(v)^* \zeta_l = \zeta_{l+m-n}$.  
So if $k_i=m_i-n_i \in K_j$, and 
$\lambda_\tau(u_i) \lambda_\tau(v_i)^* \zeta_{n_i} = \zeta_{m_i}$, then
\[
 \lambda_\tau(u_2) \lambda_\tau(v_2)^*
 \lambda_\tau(u_1) \lambda_\tau(v_1)^* \zeta_l = 
 \zeta_{l+k_1+k_2} .
\]
So $k_1+k_2 \in K_j$. Also, $-k \in K_j$ because 
$\lambda_\tau(v) \lambda_\tau(u)^*\zeta_m = \zeta_n$.

Next we note that the subgroups $K_j$ are ordered by inclusion.
That is, if $m<n$, and $\theta(\zeta_m) = \eta_i$ and $\theta(\zeta_n)=\eta_j$,
then $K_j \le K_i$.  This follows since there is a word $w\in\Fth$ so that
$\lambda_\tau(w) \zeta_m = \zeta_n$.  Consequently, $\sigma(w) \eta_i = \eta_j$.
If $k\in K_j$ and $uv^*$ is the word of degree $k$ such that 
$\sigma(u)\sigma(v)^* \eta_j = \eta_j$.
There are unique words $u',v',w'$ in $\Fth$ with the same degrees
as $u,v,w$ so that 
\[
 \eta_j = \sigma(u)\sigma(v)^*\sigma(w) \eta_i 
 = \sigma(w')\sigma(u')\sigma(v')^* \eta_i .
\]
It follows that $\sigma(u')\sigma(v')^* \eta_i$ is the unique basis vector obtained
by pulling back from $\eta_j$ by $\deg(w')=\deg(w)$ steps.
Hence $w'=w$ and $\sigma(u')\sigma(v')^* \eta_i = \eta_i$.
So $\deg(u')-\deg(v') = \deg(u)-\deg(v) = k$ belongs to $K_i$.

An increasing sequence of subgroups of $\bZ^\Bk$ is eventually constant.
So there is a subgroup $K \le H$ so that $K_j=K$ for all $\eta_j$
in a coinvariant subspace $\L$ which is the image under $\theta$ of 
$\{\zeta_m : m \in G_T\}$ for some $T \in \bZ^\Bk$.
Consider the representation $\kappa$ on $\ltwo((G_T+K)/K)$ induced by $\lambda_\tau$
obtained by collapsing cosets of $K$.  
The induced map $\tilde\theta$ of \mbox{$(G_T+K)/K$} onto the
basis of $\L$ is injective, and yields a unitary equivalence.
By Corollary~\ref{C:dilate with symmetry}, there is a group construction
representation $\rho$ on $\ltwo(\bZ^k/K)$ with symmetry group $H/K$
that dilates $\kappa$.  The minimal $*$-dilation of $\rho$ is unitarily
equivalent to the minimal $*$-dilation of $P_\L\sigma|_\L$, namely $\sigma$.
\end{proof}

Now we consider irreducibility.
Theorem~\ref{T:atomic from group} shows that it suffices to consider
group construction representations.  Thus the result we want follows from 
\cite[Theorem~6.5]{DPY} where it is established for the $\Bk=2$ case,
but the proof is not dependent on that restriction.
A group construction on $\ltwo(G)$ with symmetry group $H \ne \{0\}$
decomposes as a direct integral or sum over the dual group $\hat H$
of a family of group constructions on $\ltwo(G/H)$ with identical functions
$\ft^i$ but with different constants parameterized by $\hat H$.

\begin{thm}
An atomic $*$-representation $\sigma$ with connected graph 
is irreducible if and only if its symmetry group is trivial.

In general, if $\sigma$ is the $*$-dilation of a group construction
on $\ltwo(G)$ with symmetry group $H$, then $\sigma$ decomposes 
as a direct sum or direct integral over $\hat H$ of the $*$-dilations
of irreducible group constructions on $\ltwo(G/H)$.
\end{thm}

The import is that a complete set of the irreducible atomic $*$-repre\-sent\-ations
can be obtained by taking the inductive representation for each infinite word $\tau$,
determining the symmetry group $H$, and using this to construct a family of
irreducible group construction representations on $\ltwo(\bZ^\Bk/H)$ with 
different constants indexed by $\hat H$.

\section{Finitely Correlated Atomic Representations}\label{S:FCor}

A representation is \textit{finitely correlated} if it has a finite dimensional
coinvariant cyclic subspace.  By the previous section, this representation
must decompose as a direct sum of group constructions---and these will 
necessarily be finite groups in this case.

As in the 2-graph case, the group constructions on finite groups
are particularly tractable.  Moreover, it is possible to provide a
simple condition that determines the possible group constructions
on the product groups $G = \C_{n_1} \times \dots \times \C_{n_\Bk}$.
Every finite abelian group with $\Bk$ generators is a quotient of
$\C_{n_1} \times \dots \times \C_{n_\Bk}$, where $n_i$ is the order of $\fg_i$.

The first lemma depends only on two families of generators at a time.
So it is immediate from \cite[Lemma~7.1]{DPY}.

\begin{lem} \label{L:commute}
Let $\sigma$ be a representation of $\Fth$.
Suppose that there are words $e^i_u$ and $e^j_v$ and a unit vector $\xi$
such that 
\[
 \sigma(e^i_u)\xi=\alpha\xi  \qand   \sigma(e^j_v)\xi=\beta\xi
\]
for some $|\alpha|=|\beta|=1$.
Then $e^i_u e^j_v =e^j_v e^i_u$.
\end{lem}

\begin{cor}
If $\sigma$ is a group construction representation of $\Fth$ on a finite group
$G= \C_{n_1} \times \dots \times \C_{n_\Bk}/H$, then there are unique words
$u_i \in \Bm_i^*$ with $|u_i|=n_i$ such that $\sigma(e^i_{u_i})\xi_0 = \alpha_i\xi_0$
where $|\alpha_i|=1$ for $1 \le i \le \Bk$; and the $\{e^i_{u_i}\}$ all commute.
\end{cor}

Just as in the 2-graph case \cite{DPY}, the converse is also true.

\begin{thm}\label{k.f.c.rep}
Let $e^i_{u_i}$ be commuting words of lengths $|u_i|=n_i$,
and let $\alpha_i \in \bT$ for $1 \le i \le \Bk$. 
Then there is a unique group construction representation 
$\sigma$ of $\Fth$ on $G = \C_{n_1} \times \dots \times \C_{n_\Bk}$ 
such that $\sigma(e^i_{u_i}) \xi_0 \in \bC \xi_0$ for $1 \le i \le \Bk$
and $\alpha^i_g = \alpha_i$ for all $g \in G$.
\end{thm}

\begin{proof}
It is clear that the constants $\alpha_i$ pose no additional complication, 
so we will ignore them and assume that $\alpha_i=1$ for all $i$.

As in \cite[Lemma~7.2]{DPY}, the commutation relations of the $\Bk$-graph
completely determine the functions $\ft^i$ on $G$.
We will sketch the ideas.
Write 
\[
 u_i = j_0 j_{n_i-1} \dots j_2j_1 =: 
 \ft^i_0 \ft^i_{(n_i-1)\fg_i} \dots \ft^i_{2\fg_i} \ft^i_{\fg_i} .
\]
Then from the 2-graph case, one uses the commutation relations to
factor $e^j_{u_j} e^i_{u_i}$ in the form $e^j_v e^i_{u_{is}} e^j_w$,
where $|v|=n_j-s$ and $|w|=s$.  Consider the first factorization as two loops,
first moving in the $\fg_i$ direction starting at $\xi_0$ and returning to $\xi_0$, 
followed by the loop in the $\fg_j$ direction around to $\xi_0$ again. 
The second factorization moves first in the $\fg_j$ direction  to $\xi_{s\fg_j}$, 
then around a loop in the $\fg_i$ direction  through vectors 
$\xi_{s\fg_j + t \fg_i}$ back to $\xi_{s\fg_j}$ again, and then continuing on 
in the $\fg_j$ direction through the $\xi_{s'\fg_j}$ to $\xi_0$.  
In this way we obtain the functions $\ft^i_{s\fg_j + t \fg_i}$.  

One by one, introduce the next term $e^k_{u_k}$ and partially commute
through in order to define $\ft^i_g$ for all $g \in G$.  The fact that $\Fth$
is a $\Bk$-graph means that there is unique factorization, and so the result is
independent of the order in which this calculation is performed.
The fact that the words commute is exactly what is required in order
that one returns to the original word $e^i_{u_i}$ when one loops around,
passing this word past $e^j_{u_j}$.  
\end{proof}

\begin{example}
Choose the 3-graph in Example \ref{ffflip}. Let
$u=112$.
Then $e_u$, $f_u$, $g_u$ are mutually commuting.
Using the construction of the proof in Theorem \ref{k.f.c.rep},
we obtain a 27-dimensional finitely correlated representation $\sigma$ of
$\Fth$. It is not hard to see that $\sigma$ has a nontrivial symmetry group.
Moreover, $\sigma$ can be decomposed into a direct sum of three-dimensional atomic
representations of the following form:
\begin{alignat*}{3}
\rho(e_1)\xi_1= \omega_1 \xi_2,\ \rho(e_1)\xi_2=\omega_1 \xi_3,\ \rho(e_2)\xi_3=\omega_1 \xi_1; \\
\rho(f_1)\xi_1=\omega_2 \xi_3,\ \rho(f_1)\xi_3=\omega_2 \xi_2,\ \rho(f_2)\xi_2=\omega_2 \xi_1; \\
\rho(g_1)\xi_1=\omega_3 \xi_3,\ \rho(g_1)\xi_3=\omega_3 \xi_2,\ \rho(g_2)\xi_2=\omega_3 \xi_1.
\end{alignat*}
where $\omega_i$ are cube roots of unity.
\end{example}

We now show that there are many finitely correlated representations.

\begin{thm} \label{T:long words}
There are irreducible finite dimensional defect free representations of $\Fth$
of arbitrarily large dimension.
\end{thm}

\begin{proof}
We begin with arbitrary words $u_i \in \Bm_i^*$ and consider
$a_0 = e^1_{u_1}$ and $b_0 = e^2_{u_2} e^3_{u_3} \dots e^\Bk_{u_\Bk}$.
The technique from \cite[Proposition~7.7]{DPY} produces a pair of commuting words
as follows.  Consider the 2-graph consisting of families 
\[ E=\{e^1_u : |u|=|u_1|\} \AND F=\{ w \in \Fth : \deg(w) = (0,|u_2|,\dots,|u_\Bk| \} .\]
There is a permutation $\tilde\theta$ of $E \times F$ which determines the
commutation relation $e^1_u w = w' e^1_{u'}$ via $\tilde\theta(e^1_u,w)=(e^1_{u'},w')$.
There is a cycle beginning with $(a_0,b_0)$, namely
\[ (a_0,b_0), (a_1,b_1), \dots, (a_{n-1},b_{n-1}), (a_n,b_n) = (a_0,b_0) .\]
Then $a_ib_i = b_{i+1}a_{i+1}$ for $i \in \bZ/n\bZ$. Hence
$a := a_{n-1} \dots a_1a_0$ commutes with $b :=b_0b_1\dots b_{n-1}$.

Now we can factor $b = dc$ where $\deg(c) = (0,n|u_2|,0,\dots,0)$ 
and $\deg(d) = (0,0,n|u_3|,\dots,n|u_\Bk|)$.
At this point, we can move to the inductive step.  Suppose that we have
commuting words $a^i_0 = e^i_{v_i}$ for $1 \le i \le s$ and $c_0d_0$ where 
$c_0 = e^{s+1}_{v_{s+1}}$ and $d_0$ is a word in the remaining variables.
Consider the rank $s+2$-graph $\bF^+$ with generators 
\[ E_i = \{ e^i_u : |u|=|v_i| \}, 1 \le i \le s+1, \AND F = \{w : \deg(w) = \deg(d_0) \} .\]
We wish to build a defect free representation $\sigma$ of $\bF^+$ on $\ltwo(\bZ)$
so that for each basis vector $\delta_k$, there are unique elements $a^i_k \in E_i$
so that 
\[ \sigma(a^i_k) \delta_k = \delta_k \qfor  1 \le i \le s\AND k \in \bZ  \] 
and words $c_k \in E_{s+1}$ and $d_k \in F$ so that 
\[ \sigma(c_k) \delta_k = \delta_{k+1} \qand \sigma(d_k) \delta_{k+1} = \delta_k \qfor k \in \bZ .\]

It is clear that for such a representation to exist, certain commutation
relations must hold.
Starting at $k=0$, we must define $c_k$ and $d_k$ by the rules
\[ c_kd_k = d_{k+1}c_{k+1} \qfor k \in \bZ \]
and define $a^i_{k+1}$ by the rules
\[  c_k a^i_k = a^i_{k+1} c_k \qfor k \in \bZ. \]
We must also verify that $\{a^i_k : 1 \le i \le s \}$ commute, and
\[ d_k a^i_{k+1} = a^i_k d_k \qfor k \in \bZ . \]

To see that this does follow, observe that since $a^i_0$ commutes with $d_0c_0$,
we have $c_0a^i_0 = a^i_1 c'$; and hence
\[ d_0 a^i_1c' = d_0 c_0 a^i_0 = a^i_0 d_0 c_0 = d' a' c_0 .\]
Thus by unique factorization, $d'=d_0$, $c'=c_0$ and $a'=a^i_1$.
That is, $c_0a^i_0 = a^i_1 c_0$ and $a^i_0 d_0 = d_0 a^i_1$.
Also
\[ a^i_1 (c_0d_0) = c_0 a^i_0 d_0 = (c_0d_0) a^i_1 \]
and
\[ a^i_1 a^j_1 c_0 = c_0 a^i_0 a^j_0 = c_0 a^j_0 a^i_0 = a^j_0 a^i_0 c_0 .\]
So $\{ a^i_1 : 1 \le i \le s \} \cup \{ d_1 c_1 = c_0d_0 \}$ is a commuting family.
The relations now follow by recursion.

As before, we see that the map taking $(a^i_k, c_k,d_k)$ to $(a^i_{k+1}, c_{k+1},d_{k+1})$
results from an application of a permutation of $E_1 \times \dots \times E_{s+1} \times F$.
Hence there is an integer $n$ so that $(a^i_n, c_n, d_n) = (a^i_0, c_0, d_0)$.
Therefore, identifying $\delta_{n+j}$ with $\delta_j$, we can wrap this sequence into a
finite dimensional representation on $\spn\{\delta_k : 0 \le k < n \}$.
In particular, the words which fix $\delta_0$ must commute by Lemma~\ref{L:commute}.
So $\{ e^i_{v_i} : 1 \le i \le s \}$ together with 
$e^{s+1}_{v_{s+1}} := c_{n-1} c_{n-2} \dots c_1 c_0$ and 
$d = d_0 d_1 \dots d_{n-1}$ form a commuting family.

Repeated application of this technique produces $\Bk$ commuting words.
By Theorem~\ref{k.f.c.rep}, this gives rise to a defect free finitely correlated
representation on $\ltwo(G)$ where $G$ is a finite product group.
This may have non-trivial symmetry.
The first word had the form $e^1_{u'u_1}$ where $u_1$ was arbitrary.
We can ensure that this word has no small periods; for example, take any
word of length $N$ starting with a 1 followed by $N$ 2's.  Then this
cannot be the power of a word of length less than $2N$.
So even once we have quotiented out by the symmetry group to
obtain an irreducible representation, we have dimension at least $2N$.
\end{proof}

\section{Periodicity} \label{S:period}

In this section, we examine the periodic case in more detail.
This builds on the detailed analysis of periodicity in 2-graphs in \cite{DYperiod}.

In the case of 2-graphs, periodicity is a very special property that
requires rather stringent structural properties.  In particular,
\cite[Theorem~3.1]{DYperiod} shows that if a 2-graph $\Fth$ is periodic,
then $H_\theta = \bZ(a,-b)$ for some $a,b > 0$.
Moreover,  this occurs if and only if 
there is a bijection $\gamma: \Bm_1^a \to \Bm_2^b$ so that
\[
 e^1_u e^2_v = e^2_{\gamma(u)} e^1_{\gamma^{-1}(v)}
 \qforal u \in \Bm_1^a \AND v \in \Bm_2^b .
\]
Equivalently, the 2-graph with generators $E_1=\{e^1_u : u \in \Bm_1^a\}$
and $E_2=\{e^2_v : v \in \Bm_2^b \}$ is just a flip algebra.

The following theorem is the appropriate generalization of this result
for higher rank graphs.

\begin{thm} \label{T:period} 
Let $\Fth$ be a $\Bk$-graph, and let $a_i,b_j$ be positive integers
and $c_k = 0$ for $1 \le i \le p < j \le p+q < k \le \Bk$.
Then the following conditions are equivalent for 
$\pi=:(a_1,\dots,a_p,-b_{p+1},\dots,-b_{p+q},0,\dots,0)$
\begin{enumerate}
\item  $\Fth$ is $\pi$-periodic.
\item  Every tail is $\pi$-periodic.
\item  There is a bijection $\gamma: E \to F$, where
\[ \strut\qquad\quad
 E = \big\{ \prod_{i=1}^p e^i_{u_i} : u_i \in \Bm_i^{a_i} \big\}
 \AND
 F = \big\{ \prod_{j=p+1}^{p+q} e^j_{v_j} : v_j \in \Bm_j^{b_j} \big\},
\]
such that 
\[ef = \gamma(e) \gamma^{-1}(f) \qforal e \in E \AND f \in F ; \tag{$\dagger$} \label{dagger}  \]
and 
\[ e \tau = \gamma(e) \tau \quad\text{for every infinite tail } \tau \AND e \in E .\] 
Moreover, if $p+q=\Bk$, this condition on the tails is automatic.
\end{enumerate}
\end{thm}

\begin{proof}
If $\Fth$ is $\pi$-periodic, then by definition, every infinite tail $\tau$ is
eventually $\pi$-periodic.  Suppose that there were an infinite tail $\tau$ which 
is not $\pi$-periodic.  Then some initial segment of $\tau$, say $\tau_1$,
fails to be $\pi$-periodic.  Therefore any infinite tail which contains the
sequence $\tau_1$ infinitely often is not eventually $\pi$-periodic, contrary
to hypothesis.  Hence (ii) holds. Clearly (ii) implies (i).

Suppose that (ii) holds.  
The 2-graph $\fG^+$ with generators from $E$ and $F$
has the property that every infinite tail is $(1,-1)$-periodic.
By \cite[Theorem~3.1]{DYperiod}, 
there is a bijection $\gamma$ of $E$ onto $F$ such that 
$ef = \gamma(e) \gamma^{-1}(f)$ for all $e \in E$ and $f \in F$.
If $\tau$ is any infinite tail, take any $e\in E$ and $f\in F$, and consider the tail 
$\tau' = fe \tau = \gamma^{-1}(f)\gamma(e) \tau$. 
Since $\tau'$ is $\pi$-periodic, one obtains identical tails by deleting 
an initial word of degree $(a_1,\dots,a_p,0,\dots,0)$ or an initial word
of degree $(0,\dots,0,b_{p+1},\dots,b_{p+q},0,\dots,0)$.  That is
$e \tau = \gamma(e) \tau$.
When $p+q=\Bk$, the $\pi$-periodicity of $\tau$ is equivalent to the
$(1,-1)$-periodicity of an infinite tail in $\fG^+$; and this follows from 
the existence of $\gamma$ by \cite[Theorem~3.1]{DYperiod}.

Finally, suppose that (iii) holds. We may factor an arbitrary infinite word $\tau$
as $\tau = fe\tau'$ for some $e\in E$ and $f\in F$.
So $\tau = \gamma^{-1}(f) \gamma(e) \tau'$.
To check $\pi$-periodicity of $\tau$, we need to compare
$\gamma(e)\tau$ with $e\tau$.  These are equal by (iii). So (ii) holds.
\end{proof}

Periodicity yields non-trivial elements in the centre of $\ca(\Fth)$
just as in the 2-graph case \cite[Lemma~5.4]{DYperiod}.

\begin{cor} \label{C:centre}
Let $\Fth$ be a $\Bk$-graph with $\pi$-periodicity.  Using the notation
from Theorem~$\ref{T:period}$, define $W = \sum_{e \in E} \gamma(e) e^*$
in $\ca(\Fth)$.  Then $W$ is a unitary in the centre of $\ca(\Fth)$ satisfying
$We = \gamma(e)$ for all $e \in E$; and $W$ is a sum of terms of degree $-\pi$.
Also $W=e^*\gamma(e)$ for any $e \in E$.
\end{cor}

\begin{proof}
Since any inductive representation $\sigma$ is faithful by
Theorem~\ref{T:faithful}, we can  compute within $\B(\H_\sigma)$.
Take any standard basis vector $\xi$ and consider the infinite tail
$\tau$ obtained by pulling back from $\xi$; i.e. $\tau$ is the
unique infinite tail such that $\xi$ is in the range of $\sigma(w)$
whenever $\tau = w \tau'$.
Factor $\tau = e\tau'$ for some $e \in E$.
Let $\zeta = \sigma(e)^*\xi$ and $\xi' = \sigma(\gamma(e)) \zeta$.
Then $\xi' = \sigma(\gamma(e))\sigma(e)^* \xi = \sigma(W) \xi$.
The infinite tail obtained by pulling back from $\xi'$ is evidently
$\gamma(e) \tau'$, which equals $\tau$ by Theorem~\ref{T:period}(iii).

Now let $g= e^i_j$ be any generator of $\Fth$.  
We will show that $\sigma(W)$ commutes with $\sigma(g)^*$.
Note that $\xi$ is in the range of $\sigma(g)$ if and only if 
$\tau$ factors as $g\tau''$ (when one uses the commutation
relations to move the first $e^i_{j'}$ term to the initial position).
If $\xi$ is not in the range of $\sigma(g)$, then neither is $\xi'$,
and so 
\[ \sigma(W) \sigma(g)^*\xi = 0 = \sigma(g)^* \xi' = \sigma(g)^* \sigma(W) \xi .\]
If $\xi$ is in the range of $\sigma(g)$, then there is a unique word $e' \in E$
so that $\tau = g e' \tilde\tau$.  Also factor $ge' = eg'$ for some $g'=e^i_{j'}$.
Let $\eta = \sigma(g)^* \xi$, $\eta' = \sigma(g)^* \xi'$
and $\zeta' = \sigma(e')^* \eta = \sigma(g')^* \zeta$.
Now $\tau = g \gamma(e') \tilde\tau$  by Theorem~\ref{T:period}(iii).
Thus $\xi'$ is in the range of $\sigma(g \gamma(e'))$, say 
$\xi' = \sigma(g \gamma(e')) \tilde{\zeta}$.  Now
$\gamma(e) g' = \tilde g \gamma(\tilde e)$ for some $\tilde g = e^i_{\tilde j}$
and $\tilde e \in E$.  But then
\[ 
 \sigma(g \gamma(e')) \tilde\zeta = \xi' = 
 \sigma(\gamma(e) g') \zeta' = \sigma(\tilde g \gamma(\tilde e)) \zeta' .
\]
It follows that $\tilde g = g$, $\tilde e = e'$ and $\tilde{\zeta} = \zeta'$.
Therefore 
\[
 \sigma(g)^* \sigma(W) \xi = \sigma(g)^* \xi' = \eta' = \sigma(\gamma(e'))\sigma(e')^* \eta
 = \sigma(W) \sigma(g)^* \xi .  
\]
We conclude that $W$ commutes with $g^*$ for every generator of $\Fth$.

It is evident that $\sigma(W)$ carries the range of each 
$\sigma(e)$ for $e\in E$ onto the range of $\sigma(\gamma(e))$.
As the ranges of $\sigma(e)$ for $e\in E$ are pairwise orthogonal 
and sum to the whole space, as do the ranges of $\sigma(\gamma(e))$,
it follows that $\sigma(W)$ is unitary.  
Thus it commutes with all of $\ca(\Fth)$; i.e. it is in the centre.
The identity $We = \gamma(e)$ for $e \in E$ is clear.
Also $\deg(\gamma(e) e^*) = \deg(\gamma(e)-\deg(e) = -\pi$ by construction.

Now $eW=We=\gamma(e)$, and hence $W=e^*\gamma(e)$.
\end{proof}

We now take a closer look at this for $\Bk = 3$.
The permutation $\theta_{12}$ which define the relations between
$\{e^1_i\}$ and $\{e^2_j\}$ determines a function $\tilde\theta_{12}$
from $\Bm_1^*\times \Bm_2^*$ so that $\tilde\theta_{12}(u,v) = (u',v')$
when $e^1_ue^2_v = e^2_{v'}e^1_{u'}$.

\begin{prop} \label{3Gperiodic}
Let $\Fth$ be a  $3$-graph with generators $\{e_i : i \in \Bl\}$,
\mbox{$\{f_j : j \in \Bm\}$} and $\{g_k : k \in \Bn\}$.
Suppose that $\Fth$ has $(a,b,-c)$ symmetry,
where $a,b,c \in \bN$; and let  $\gamma : \Bl^a \times \Bm^b \to \Bn^c$
satisfy {\em ($\dagger$)}, i.e.
\[
 e_{u_0} f_{v_0} g_{\gamma(u_1,v_1)} = g_{\gamma(u_0,v_0)} e_{u_1}f_{v_1} 
 \qforal (u_i,v_i) \in \Bl^a \times \Bm^b .  
\]
Let $\delta = \gamma \tilde\theta_{12}^{-1}$. 
Then for all $(u_i,v_i) \in \Bl^a \times \Bm^b$,
\[
 e_{u_0} g_{\delta(u_1,v_0)} = g_{\gamma(u_0,v_0)}e_{u_1}  
 \qand
 f_{v_0} g_{\gamma(u_1,v_1)} = g_{\delta(u_1,v_0)}f_{v_1} ;
\]
and conversely these relations imply that $\Fth$ is $(a,b,-c)$-periodic.
Moreover,
\[  f_{v_0} e_{u_0} g_{\delta(u_1,v_1)} = g_{\delta(u_0,v_0)}f_{v_1}e_{u_1} .\]
\end{prop}

\begin{proof}
There is a word  $w \in \Bn^c$ so that
\[
 g_{\gamma(u_0,v_0)} e_{u_1}f_{v_1} = e_{u_0}g_{w} f_{v_1} = 
 e_{u_0} f_{v_0}g_{\gamma(u_1,v_1)} .
\]
Since $g_{\gamma(u_0,v_0)} e_{u_1} = e_{u_0}g_{w}$, 
$w$ depends only on $u_0,u_1,v_0$;
and since $g_{w} f_{v_1} = f_{v_0}g_{\gamma(u_1,v_1)}$, 
$w$ depends only on $u_1,v_0, v_1$.
Therefore $w$ is a function $w = \delta(u_1,v_0)$.  That is,
\[
 e_{u_0} g_{\delta(u_1,v_0)} = g_{\gamma(u_0,v_0)}e_{u_1}  
 \qand
 f_{v_0} g_{\gamma(u_1,v_1)} = g_{\delta(u_1,v_0)}f_{v_1} .
\]
Hence
\[
 f_{v_0} e_{u_0} g_{\delta(u_1,v_1)} = f_{v_0} g_{\gamma(u_0,v_1)} e_{u_1}
 = g_{\delta(u_0,v_0)} f_{v_1}  e_{u_1} . 
\]

Fix $(u_0,v_0)$ and $(u_1,v_1)$.  
Let $\tilde\theta_{12}^{-1}(u_i,v_i) = (u_i',v_i')$ for $i=0,1$;
so that $f_{v_i}e_{u_i} = e_{u'_i}f_{v'_i}$. 
Also let $\delta^{-1}\gamma(u'_1,v'_1) = (u,v)$. Then
\begin{align*}
 g_{\gamma \tilde\theta_{12}^{-1}(u_0,v_0)} f_{v_1}e_{u_1} &=
 g_{\gamma(u'_0,v'_0)} e_{u'_1} f_{v'_1} = 
 e_{u'_0}f_{v'_0} g_{\gamma(u'_1,v'_1)} \\ &=
 f_{v_0} e_{u_0} g_{\delta(u,v)} =
 f_{v_0} g_{\gamma(u_0,v)} e_u \\ &=
 g_{\delta(u_0,v_0)} f_v e_u .
\end{align*}
Hence $\delta(u_0,v_0) = \gamma \tilde\theta_{12}^{-1}(u_0,v_0)$.

The converse is straightforward.
By Theorem~\ref{T:period} (iii), the tail condition is automatic.
Hence $\Fth$ is $(a,b,-c)$-periodic.
\end{proof}

Proposition \ref{3Gperiodic} allows us to define some examples of periodic 3-graph algebras.

\begin{example}
Suppose that $\theta_{12}=\id$; i.e.\ the $e_i$'s commute with the $f_j$'s.
Also suppose that $n = lm$ and fix a bijection $\gamma:\Bl\times \Bm \to \Bn$.
By Proposition \ref{3Gperiodic}, we require $\delta=\gamma$.
So define relations
\[
 e_i g_{\gamma(i',j')} = g_{\gamma(i,j')} e_{i'} \qand
 f_j g_{\gamma(i',j')}=g_{\gamma(i',j)}f_{j'}.
\]
It is easy to check the cubic condition---so this determines a $3$-graph $\Fth$.
By Proposition \ref{3Gperiodic}, $\Fth$ is $(1,1,-1)$-periodic.
Its symmetry group is exactly $\bZ(1,1,-1)$.
\end{example}

\begin{example}
Suppose that $l=m$, $n=m^2$ and $\theta_{12}$ is the transposition:
$\theta_{12}(i,j)=(j,i)$; i.e. $e_if_j = f_ie_j$. Then
we again identify a bijection
$\gamma:\Bl \times \Bm \to \Bn$ and motivated by Proposition \ref{3Gperiodic}, 
define $\delta(i,j) = \gamma(j,i)$. Then define the commutation relations
\[
 e_{i}g_{\gamma(j,k)} = g_{\gamma(i,j)}e_{k} \qand
 f_{i}g_{\gamma(j,k)} = g_{\gamma(i,j)}f_{k}.
\]
This is easily seen to be a $3$-graph $\Fth$.
By construction, it has $(1,1,-1)$-periodicity.

The first two variables form a 2-graph with $(1,-1)$-periodicity.
Since $(1,-1,0)$ has a 0, it is convenient to look for 
$(3,1,-2)$-periodicity instead.  One readily computes that
\begin{align*}
e_{i_1}e_{i_2}e_{i_3}f_j g_{\gamma(i'_1,i'_2)} g_{\gamma(i'_3,j')} &=
e_{i_1}e_{i_2} g_{\gamma(i_3,j)} e_{i'_1}f_{i'_2} g_{\gamma(i'_3,j')} \\&=
e_{i_1} g_{\gamma(i_2,i_3)} e_j g_{\gamma(i'_1,i'_2)} e_{i'_3}f_{j'} \\&=
g_{\gamma(i_1,i_2)} e_{i'_3} g_{\gamma(j,i'_1)} e_{i'_2} e_{i'_3}f_{j'} \\&=
g_{\gamma(i_1,i_2)} g_{\gamma(i_3,j)} e_{i'_1} e_{i'_2} e_{i'_3}f_{j'}
\end{align*}
Thus the function $\Gamma:\Bl^3\times\Bm \to \Bn^2$ by
$\Gamma(i_1,i_2,i_3,j) = (\gamma(i_1,i_2), \gamma(i_3,j))$
is a bijection satisfying $(\dagger)$.

Thus the symmetry group is 
\[ H = \bZ(1,1,-1) + \bZ(3,1,-2) =  \bZ(1,-1,0) + \bZ(1,1,-1) .\]
\end{example}

\begin{example}
A different example can be defined when $l=m=n=2$.
Let $\Fth$ be the 3-graph with $\bF^+_{\theta_{12}}$ being the
flip algebra, and $\bF ^+_{\theta_{13}}$, $\bF ^+_{\theta_{23}}$ 
each being the square algebra.
A straightforward computation shows that 
$\Fth$ has the symmetry group $H=\bZ(1,-1,0) +\bZ(2,0,-2)$.
\end{example}

\begin{example}
In \cite{DYperiod}, many examples of periodic 2-graphs were exhibited, 
some with surprisingly high order of periodicity.
It is easy to combine a number of 2-graphs together with other
variables by making them commute.
So suppose that $\bF^+_{\theta_{2i-1,2i}}$ are 2-graphs with symmetry groups 
$H_i = \bZ(a_i,-b_i)$ for $1 \le i \le s$, and let $\Fth$ be any $k$-graph 
with symmetry group $H \le \bZ^k$.
Form a $2s+k$-graph $\fG^+_\theta = \prod_{i=1}^s \bF^+_{\theta_{2i-1,2i}} \times \Fth$
by declaring that variables in the different factors of the product commute.
Then it is routine to check that this is a $2s+k$--graph with symmetry group
$H_\theta = \prod_{i=1}^s H_i \times H \le \bZ^{2s+k}$.
\end{example}

\section{The Structure of Graph C*-algebras}

Kumjian and Pask \cite{KumPask} showed that $\ca(\Fth)$ is simple 
if $\Fth$ is aperiodic.
Robertson and Sims \cite{RS} proved the converse.
In \cite{DYperiod}, we showed that in the case of a periodic 2-graph
$\Fth$ on a single vertex, one has the more precise description 
$\ca(\Fth) \cong \rC(\bT) \otimes \fA_1$ for some simple C*-algebra $\fA_1$.
We wish to extend this result to $\Bk$-graphs (on a single vertex).
The proof will follow the method of \cite{DYperiod} of constructing
two approximately inner expectations.

Let $H_\theta$ be the symmetry group of $\Fth$. 
Then $H_\theta \cong \bZ^s$ for some $s \le \Bk$.
Let $h_1,\dots,h_s$ be a set of free generators, and set $\vec h = (h_1,\dots,h_s)$. 
By Corollary~\ref{C:centre}, each $h \in H_\theta$ determines
a unitary operator $W_h$ in the centre $\Z$ of $\ca(\Fth)$.
It turns out that this map is a group homomorphism---but we
will not establish that at this time.  Instead we define $W_i = W_{h_i}$
and use these as generators for an abelian algebra which will
eventually turn out to be all of $\Z$.
For $n =(n_1,\dots,n_s) \in \bZ^s$, write $W^n = \prod_{i=1}^s W_i^{n_i}$ 
and $n \cdot \vec h = \sum_{i=1}^s n_ih_i$. \vspace{.3ex}

Let $G_\theta = \bZ^\Bk/H_\theta$, and let $\pi$ be the quotient map.
Recall that every character $\phi \in \wh{\bZ^\Bk} \cong \bT^\Bk$
determines a gauge automorphism $\gamma_\phi$ of $\ca(\Fth)$ such that
$\gamma_\phi(w) = \phi(\deg w) w$ for all $w \in \Fth$.
Each character $\psi \in \wh{G_\theta}$ determines the character $\psi\pi$
on $\bZ^\Bk$ which takes the value 1 on $H_\theta$.
Let $\gamma_\psi$ denote the corresponding gauge automorphism.
We will define an expectation on $\ca(\Fth)$ which respects the symmetry of $\Fth$:
\[ \Psi(X) = \int_{\wh{G_\theta}} \gamma_\psi(X) \,d\psi .\]

\begin{lem}
The joint spectrum of $(W_1,\dots,W_s)$ is $\bT^s$; and
\[ \ca(\fF, W_1,...,W_s)\cong \rC(\bT^s)\otimes\fF\cong \rC(\bT^s,\fF) . \]
\end{lem}

\begin{proof}
If $\upchi$ is any character of $H_\theta$, there is a character
$\phi$ in $\wh{\bZ^\Bk}$ which extends $\upchi$.
Therefore
\[ \gamma_\phi(W^n) = \phi(n\cdot\vec h) W^n = \upchi(n) W^n .\]
Thus $\gamma_\phi$ restricts to an automorphism of $\ca(W_1,\dots,W_s)$
and for any $\lambda_i \in \bT$, there is some $\phi$ so that
$\gamma_\phi(W_i) = \lambda_i W_i$ for $1 \le i \le s$.
So the joint spectrum $\sigma(W_1,\dots,W_s)$ is invariant under
the transitive action of the torus.
Hence $\sigma(W_1,\dots,W_s) = \bT^s$ and
$\ca(W_1,\dots,W_s) \cong \rC(\bT^s)$.

There is a canonical map of the tensor product 
$\rC(\bT^s) \otimes \fF \cong \rC(\bT^s,\fF)$
onto $\ca(\fF, W_1,...,W_s)$ which sends the constant functions onto $\fF$
and sends $z_i$ to $W_i$.  Since $\fF$ is simple, the kernel consists of
all functions vanishing on some closed subset of $\bT^s$.  However,
since the joint spectrum of $(W_1,\dots,W_s)$ is all of $\bT^s$, this set
must be empty; and this map is an isomorphism.
\end{proof}

\begin{thm} \label{T:expectation}
The map $\Psi$ is a faithful, completely positive, approximately inner expectation onto
$\ca(\fF, W_1,...,W_s)$.
\end{thm}

\begin{proof} 
Since $\Psi$ is the average of automorphisms, it is clearly faithful and
completely positive.

Let $w = uv^*$, where $u,v \in \Fth$, be a typical word in $\ca(\Fth)$.
Then $\deg(w) = \deg(u)-\deg(v)$ is a homomorphism into $\bZ^\Bk$.
The kernel, the words of zero degree, generate $\fF$ as a C*-algebra.
If $\deg(w) = h \in H_\theta$, then $\pi\deg(w) = 0$;
so for any $\psi \in \wh{G_\theta}$, 
we have $\gamma_\psi(w) = w$. Therefore $\Psi(w) = w$.
On the other hand, if $\deg(w) \not\in H_\theta$, then $\pi\deg(w) = g \ne 0$.
Therefore 
\[ \Psi(w) = \int_{\wh{G_\theta}} \psi(g) \,d\psi\ w = 0 .\]

By Corollary~\ref{C:centre},$W_i$ is a sum of words of degree $-h_i$.
Hence $W^n$ is a sum of words of degree $-n\cdot \vec h$.
Thus if $\deg(w) = h = n \cdot \vec h \in H_\theta$, 
then $w W^n$ is a sum of words of degree 0, and so lies in $\fF$.
Therefore $w = F W^{*n}$ belongs to $\ca(\fF, W_1,...,W_s)$.
Conversely, this C*-algebra is spanned by terms of the form
$uv^* W^n$ where $\deg(uv^*)=0$ and $n \in \bZ^s$. We have
seen that $\Psi(uv^* W^n) = uv^* W^n$.  Thus $\Psi$ is an expectation
onto $\ca(\fF, W_1,...,W_s)$.

Lastly, we will construct a sequence of isometries $V_n \in \ca(\Fth)$
so that $\lim_{n\to\infty} V_n^* X V_n = \Psi(X)$ for all $X \in \ca(\Fth)$.

By Proposition~\ref{P:goodtail}, there is an infinite tail $\tau$ with
$H_\tau = H_\theta$.
Therefore if $u,v \in \Fth$ and $\deg(uv^*) \not\in H_\theta$, then
$u\tau \ne v\tau$.  It follows that for some initial segment $\tau_0$
of $\tau$, $\tau_0^* u^*v \tau_0 = 0$.
So for each integer $n \ge 1$, select $\tau_n$ so that
$\tau_n^* (u^*v) \tau_n = 0$ whenever $\deg(uv^*) \not\in H_\theta$
and \mbox{$\deg(u) \vee \deg(v) \le \Bn := (n,n,\dots,n)$.}

Let $\S_n = \{ x \in \Fth : \deg(x) = \Bn \}$.
Define an isometry in $\ca(\Fth)$ by 
\[ V_n = \sum_{x \in \S_n} x \tau_n x^* .\]

Suppose that $uv^* \in \fF$; so $\deg(u)=\deg(v)$.
For $n$ sufficiently large, $\deg(u)  \le \Bn$.
Then $uv^*$ can be written as a sum of words of the same form
with $\deg(u)=\deg(v) = \Bn$.
Recall that when $\deg(x)=\deg(u)$, then $x^*u = \delta_{x,u}$.
Therefore
\begin{align*}
 V_n^* (uv^*) V_n &= \sum_{x \in \S_n} \sum_{y \in \S_n}
 x \tau_n^* (x^* u) (v^* y) \tau_n y^*  =
 u \tau_n^* \tau_n v^* = uv^* .
\end{align*} 
It follows that $V_n^*(uv^*)V_n = uv^*$ whenever 
$\deg(u)=\deg(v) \le \Bn$.

Next suppose that $\deg(uv^*) = h = m \cdot \vec h \in H_\theta$.
Then $W^m(uv^*)$ belongs to $\Fth$.
By the previous paragraph, for $n$ sufficiently large,
\begin{align*}
 V_n^*(uv^*)V_n &= V_n^*W^{*m}(W^muv^*)V_n = 
 W^{*m} V_n^*(W^muv^*)V_n \\ &= W^{*m}(W^muv^*) = uv^* .
\end{align*} 

Finally, suppose that $\deg(uv^*) \not\in H_\theta$.
For sufficiently large $n$, we have $\deg(u) \vee \deg(v) \le \Bn$.
Hence if $x,y \in \S_n$, then $x^*(uv^*)y$ is either $0$ or
it has the form $ab^*$ where $\deg(a) \vee \deg(b) \le \Bn$
and
\[ \deg(a)-\deg(b) = \deg(u)-\deg(v) \not\in H_\theta .\]
Consequently, 
\[ x\tau_n^*x^*(uv^*)y\tau_ny^* = x(\tau_n^*ab^* \tau_n)y^* = 0 .\]
Summing over $\S_n\times\S_n$ yields $V_n^*(uv^*)V_n=0$.

We have shown that 
\[ \lim_{n\to\infty} V_n^* (uv^*) V_n = \Psi(uv^*) \]
for all words $uv^*$.  Thus this limit extends to all of $\ca(\Fth)$.
\end{proof}

The next step is to localize at the point $\one = (1,\dots,1)$ of $\bT^s$
in $\sigma(W_1,\dots,W_s)$.  Essentially this is evaluation at $\one$.

\begin{thm}\label{commutingd}
Let $q$ be the quotient map of $\ca(\Fth)$ onto
\[ \fA:=\ca(\Fth)/\langle W_1-I, ..., W_s-I\rangle .\]
Let $\ep_\one$ be evaluation at $\one$ in $\rC(\bT^s,\fF)$.
Then there is a faithful, completely positive, 
approximately inner expectation $\Psi_1$ of
$\fA$ onto $\fF$ such that the following diagram commutes:
\[
 \xymatrix{
\ca(\Fth) \ar[r]^(.6){q} \ar[d]_{\Psi} & \fA \ar@{->}[d]^{\Psi_1}\\
 \rC(\bT^s,\fF) \ar[r]_(.6){\ep_\one}&\fF
}
\]
Moreover, $\fA$ is a simple C*-algebra.
\end{thm}

\begin{proof}
For each $\psi \in \wh{G_\theta}$, we have $\gamma_\psi(W_i) = W_i$.
Hence $\gamma_\psi$ carries the ideal $\langle W_1-I, ..., W_s-I\rangle$ onto itself.
Therefore it induces an automorphism $\dot\gamma_\psi$ of $\fA$;
and $\dot\gamma_\psi q = q \gamma_\psi$.
Define a map on $\fA$ by
\[
 \Psi_1(A) = \int_{\widehat{G_\theta}} \dot\gamma_\psi (A) \,d\psi
 \qfor A\in\fA.
\]
Then it follows that $\Psi_1 q = q \Psi = \ep_\one \Psi$ because
the restriction of $q$ to $\ca(\fF,W_1,\dots,W_s) \cong \rC(\bT^s,\fF)$
is evidently $\ep_\one$.

Since $\Psi_1$ is an average of automorphisms,
it is a faithful, completely positive map into $\fF$.
Now $q$ is an isomorphism on $\fF$; so we may identify
$q\fF$ in $\fA$ with $\fF$. For any $F \in \fF$, 
\[ \Psi_1(F) = \ep_\one(\Psi(F)) = \ep_\one(F) = F .\] 
Thus $\Psi_1$ is an expectation.

Let $\dot V_n = q(V_n)$.
$\Psi_1$ is approximately inner because if $A = qX \in \fA$,
\[
 \Psi_1(A) = q \Psi(X) =
 \lim_{n\to\infty} q( V_n^* X V_n) =
 \lim_{n\to\infty} \dot V_n^* A \dot V_n . 
\]

The fact that $\fA$ is simple now follows.
Indeed, if $\fJ$ is a non-zero ideal of $\fA$, then it
contains a non-zero positive element $A$.
Since $\Psi_1$ is faithful, $\Psi_1(A) \ne 0$.
Also because $\Psi_1$ is approximately inner,
$\Psi_1(A)$ belongs to $\fF \cap \fJ$.
Since $\fF$ is simple, we conclude that
$\fF \cap \fJ = \fF$ contains the identity; and
hence $\fJ = \fA$.
\end{proof}

\begin{rem}
The expectation $\Psi_1 q = \ep_\one \Psi$ of $\ca(\Fth)$
onto $\fF$ is not the expectation $\Phi$ obtained by integrating over
all gauge automorphisms.  Indeed $\Phi(W_i)=0$ while
$\Psi_1q(W_i) = \Psi_1(I) = I$.
\end{rem}

We are now ready to obtain the structure of $\ca(\Fth)$.

\begin{thm}\label{T:tensor}
Suppose that the symmetry group $H_\theta$ of $\Fth$ has \mbox{rank $s$.}
Then 
\[ \ca(\Fth) \cong \rC(\bT^s) \otimes \fA. \]
\end{thm}

\begin{proof}
For $z=(z_1,\dots,z_s) \in \bT^s$, let $\fJ_z = \langle W_1-z_1I,...,W_s-z_sI\rangle$
and set $\fA_z:=\ca(\Fth)/\fJ_z$.
Let $\upchi_z \in \wh{H_\theta}$ be the character $\upchi_z(h_i) = z_i$,
and let $\phi$ be any extension to a character on $\wh{\bZ^\Bk}$.
Since $\deg(W_i) = -h_i$, we have $\gamma_\phi(W_i) = \bar z_i W_i$.
Therefore $\gamma_\phi(W_i - I) = \bar z_i (W_i - z_i I)$.
So $\gamma_\phi$ carries the ideal $\fJ_\one$
onto $\fJ_z$.
It follows that $\fA_z \cong \fA$ via the automorphism $\dot\gamma_\phi$.

Moreover, the lifting from $\upchi$ to $\phi$ can be done in a 
locally continuous way.  That is, for any $z \in \bT^s$, there is a 
neighbourhood of $z$ on which we can select a continuous lifting.
It follows that the map taking $z$ to $q_z(X)=q(\gamma_\phi(X))$ 
is a continuous function into $\fA$.
This provides a $*$-homomorphism $\Theta$ of $\ca(\Fth)$ into $\rC(\bT^s,\fA)$.

The restriction of $\Theta$ to $\rC(\bT^s,\fF)$ is readily seen to be the
identity map.  In particular, this restriction is an isomorphism.
By Theorem~\ref{commutingd}, it follows that 
$\Theta \Psi = (\Psi_1 \otimes \id) \Theta$.
The left hand side is faithful, and therefore $\Theta$ must 
be an monomorphism.

If $A \in \fA$ and $qX=A$, then for any $f \in \rC(\bT^s)$, 
\[ \Theta( f(W_1,\dots,W_s) X) = f(z_1,\dots,z_s) A .\]
These functions span $\rC(\bT^s,\fA)$, and therefore $\Theta$
is surjective.  This establishes the desired isomorphism.
\end{proof}

As an immediate consequence, we
obtain the following characterization of the simplicity of $\Bk$-graph C*-algebras,
which was proved in full generality by Kumjian--Pask \cite{KumPask} (sufficiency) and
Robertson-Sims \cite{RS} (necessity) using completely different approaches.

\begin{cor}
$\ca(\Fth)$ is simple if and only if $\Fth$ is aperiodic.
\end{cor}

As another immediate consequence, we have identified the center of $\ca(\Fth)$.

\begin{cor}
The center of $\ca(\Fth)$ is $\ca(W_1,\dots,W_s) \cong \rC(\bT^s)$.
\end{cor}

Corollary~\ref{C:centre} defined an element $W_h$ in $\Z$ for any $h\in H_\theta$.
If $h = (a_1,\dots,a_n)$, then we set 
$E = \big\{ \prod e^i_{u_i} : a_i > 0,\ u_i \in \Bm_i^{a_i} \big\}$
 and
 $F = \big\{ \prod e^i_{v_i} : a_i < 0,\ u_i \in \Bm_i^{-a_i} \big\}$.
There is a bijection $\gamma$ of $E$ onto $F$ so that 
$ef = \gamma(e) \gamma^{-1}(f)$ for $e \in E$ and $f \in F$.
The central unitary $W_h$ is then defined
by the formula $W_h = \sum_{e \in E} \gamma(e)e^*$.
This is the sum of words of degree $-h$.

Write $h = n \cdot \vec h$.
In $\ca(W_1,\dots,W_s)$, the unitary $W^n$
also has $\deg(W^n) = -h$.
Therefore, in the identification of $\ca(W_1,\dots,W_s)$ with $\rC(\bT^s)$
which sends $W_i$ to $z_i$, both $W_h$ and $W^n$ are sent to a scalar
multiple of $z^n$.  
However, it is clear from the commutation relations in $\Fth$
that multiplying out the product $W_1^{n_1} \dots W_s^{n_s}$ will
be a sum of words of the form $uv^*$ with no scalars, as is $W_h$.
So they must be equal.

Since $\Z = \ca(W_1,\dots,W_s) = \spn\{ W^n : n \in \bZ^s \}$,
we deduce:

\begin{cor}
The map taking $h \in H_\theta$ to $W_h$ is a group homomorphism
into the unitary group of $\Z$; and $\Z = \spn\{ W_h : h \in H_\theta\}$.
\end{cor}


\end{document}